\documentclass[preprintnumbers,amsmath,amssymb]{revtex4}

\usepackage{epsfig}
\usepackage{graphicx}% Include figure files
\usepackage{dcolumn}% Align table columns on decimal point
\usepackage{bm}% bold math
\usepackage{tikz}
\usepackage[all]{xy}
\usetikzlibrary{arrows,shapes}
\usepackage{url}
\usepackage{centernot}
\usepackage{stmaryrd}

\usepackage{color} % REMOVE BEFORE SUBMITTING
\usepackage{fancyheadings} 
\fancypagestyle{first}{% 
\fancyhf{} % clear all header and footer fields 
\cfoot{
Approved for public release; unlimited distribution. \copyright 2017 BAE Systems. All rights reserved.
} % except the center 
 
} 

\newcommand{\tstamp}{\today}   
\makeatletter
\newcommand{\xMapsto}[2][]{\ext@arrow 0599{\Mapstofill@}{#1}{#2}}
\def\Mapstofill@{\arrowfill@{\Mapstochar\Relbar}\Relbar\Rightarrow}
\makeatother

\begin{document}

\title{Path abstraction}% Force line breaks with \\

\author{Steve Huntsman}%
\email{steve.huntsman@baesystems.com}
\affiliation{BAE Systems, 4301 North Fairfax Drive, Arlington, Virginia 22203
}%

\date{\today}% It is always \today, today,

\begin{abstract}
Given the set of paths through a digraph, the result of uniformly deleting some vertices and identifying others along each path is coherent in such a way as to yield the set of paths through another digraph, called a \emph{path abstraction} of the original digraph. The construction of path abstractions is detailed and relevant basic results are established; generalizations are also discussed. Connections with random digraphs are also illustrated.
\end{abstract}


\maketitle

%% HEADER/FOOTER FIRST PAGE
\thispagestyle{first}

%keywords: yadda, blah, whatever

\section{\label{sec:Introduction}Introduction}

Each path in a digraph $D$ corresponds to a word over the alphabet $V(D)$. Given a subset $U \subseteq V(D)$, consider the map $\nabla_U$ that deletes elements of $U$ from such words. The question naturally arises: does the image of $\nabla_U$ correspond to the set of paths in some other digraph? 

In this paper, we address this and related questions. The practical motivation is that we have a complicated structure such as a digraph representing the possible flow of some some quantity through a system, and we wish to abstract away irrelevant details while efficiently preserving paths in the structure. For example, we might consider the flow of data \cite{NielsonNielsonHankin}, taint \cite{SchwartzAvgerinosBrumley} or provenance \cite{CheneyAcarPerera} in computer programs or systems. In general, we cannot assume that the structure has a hierarchical or modular organization, and so clustering or decomposition techniques do not solve the task at hand. Instead, we introduce a natural construction (originally proposed by Mukesh Dalal) that can be used in many circumstances to delete irrelevant vertices and identify related vertices (and though of course clustering and decomposition techniques can have something to say in the determination of these vertices, this issue will not be considered here). Subsequently, we can reason about paths in this construction rather than in its larger antecedent structure.

The paper is organized as follows: in \S \ref{sec:Basics}, we introduce basic notation and definitions; \S \ref{sec:VertexAbstraction} is a sort of appetizer from the point of view of vertices rather than paths; \S \ref{sec:DetoursBypassesAbstractions} introduces the key constructions for digraphs, which are generalized in \S \ref{sec:WeightedDigraphs} to weighted digraphs. Random digraphs are considered in \S \ref{sec:RandomDigraphs}. So-called temporal networks that are essentially time series of arcs are considered in \S \ref{sec:TemporalNetworks}. Finally, appendices give alternative proofs of key results for digraphs and briefly indicate the potential relevance of our constructions to renormalization and percolation on random digraphs.

\section{\label{sec:Basics}Basic notation and definitions}

In this paper we generally follow (or at least adapt in an obvious way) the definitions and notation of \cite{BangJensenGutin} without further comment. In particular, digraphs are loopless, so digraph morphisms are unambiguously defined in terms of their action on vertices. Also, $A$ indicates a set of arcs; the adjacency matrix corresponding to a colored/labeled digraph, directed multigraph, or weighted digraph $D$ is $\mu_D$. We assume that $\mu_D$ takes values in an appropriate commutative semiring, e.g., the Boolean semiring for digraphs. We use $+$ and $\cdot$ to denote both ordinary arithmetic and generic semiring operations, while $\lor$ denotes either logical disjunction or maximum depending on context; similarly, $\land$ denotes either logical conjunction or minimum.
\footnote{
Examples of semirings are the Boolean semiring $(\{0,1\}, \lor, \land)$ for digraphs and the real semiring $(\mathbb{R}, +, \cdot)$ for weighted digraphs. Other examples of potential relevance are the tropical semiring $(\mathbb{R} \cup \infty, \wedge, +)$ \cite{SpeyerSturmfels}, the fuzzy semiring $([0,1], \vee, \wedge)$, the bottleneck semiring $(\mathbb{R}, \vee ,\wedge)$, the log semiring $(\mathbb{R} \cup \infty, +_{\log}, +)$ with $x +_{\log} y := -\log(\exp(-x)+\exp(-y))$, the Viterbi semiring $([0,1], \vee, \cdot)$, etc. For general background on semirings, see \cite{Glazek}. There is no truly standard definition of a semiring, and the same appears to be true of related notions. Rather than invoking necessarily idiosyncratic terminology to remove any ambiguity here, we (perhaps atypically) choose to rely on common sense. For instance, we will require a zero element to indicate the absence of arcs, but this need not be a neutral element for the additive operation. Thus, e.g., we can consider the restriction of the usual tropical semiring to $[0,\infty)$. 
}

For a digraph $D = (V,A)$ and vertex coloring or labeling $\ell : V \rightarrow \Lambda$, write $D(\ell) := (\ell,A)$ for the corresponding colored digraph, omitting the dependence on $\ell$ if the desire exists and context allows. Without loss of generality, we shall assume that $V(D) = [n] \equiv \{1,\dots,n\}$ for some $n \in \mathbb{N}$.

As a shorthand, for $U \subseteq V(D)$, we shall write $D / U \equiv D / D \langle U \rangle$ for the vertex contraction of $U$ in $D$. That is, we set $V(D / U) := (V(D) \setminus U) \cup \{U\}$ and for $x,y \in V(D) \setminus U$,
\begin{eqnarray}
\mu_{D / U}(x,\{U\}) & := & \sum_{u \in U} \mu_D(x,u), \nonumber \\
\mu_{D / U}(\{U\},x) & := & \sum_{u \in U} \mu_D(u,x), \nonumber \\
\mu_{D / U}(x,y) & := & \mu_D(x,y). \nonumber
\end{eqnarray}
Similarly, for disjoint subsets $\{U_j\}_{j \in [m]}$ of $V(D)$, $D / \{U_1,\dots,U_m\} := (\dots (D / U_1)\dots) / U_m$ is well-defined.

\section{\label{sec:VertexAbstraction}Vertex abstraction}

A subset $L = \{L_j\}_{j \in [m]} \subseteq \Lambda$ of colors determines a \emph{partial partition} $\pi_{(\ell,L)}$ of $[n]$ (i.e., a partition of a subset of $[n]$) as follows: the blocks $\pi_{(\ell,L)}^{(j)}$ of $\pi_{(\ell,L)}$ are simply the preimages $\ell^{-1}(L_j)$ for $j \in [m]$. Conversely, any partial partition of $[n]$ is nonuniquely determined by some $\pi_{(\ell,L)}$, but we can choose a canonical representative for $\ell$ that makes the correspondence between colorings and partial partitions bijective.
\footnote{
Given $\ell$, consider $[\ell] := \{ f : \text{dom } f = [n] \wedge \pi_{(f,f[n])} = \pi_{(\ell,\Lambda)} \}$. By construction, $[\ell]$ is the equivalence class of functions yielding the same partition of $[n]$. We can construct a canonical representative $\ell'$ of $[\ell]$ as follows: let $\ell' : [n] \rightarrow [|\ell([n])|]$ be such that $\min \ell'^{-1}(j) < \min \ell'^{-1}(j+1)$ for $j \in [|\ell([n])|-1]$. For example, if (using a standard notation for the Cartesian product of functions for brevity) $\ell^{\times 5}(1,2,3,4,5) = (4,2,0,2,0)$, the corresponding canonical representative $\ell'$ is defined by $\ell'^{\times 5}(1,2,3,4,5) = (1,2,3,2,3)$. This canonical representative is essentially the so-called \emph{restricted growth string} corresponding to $\pi_{(\ell,\Lambda)}$ \cite{Knuth}.
} 
Henceforth we shall assume without loss of generality that $\ell$ is canonical (and so also surjective) unless otherwise specified.

As usual, let $\Pi_n$ denote the lattice of partitions of $[n]$. Following \cite{HanlonHershShareshian}, we consider the lattice $\Pi_{\le n}$ of partial partitions ordered by refinement, i.e., for $\pi, \pi' \in \Pi_{\le n}$, we have $\pi \le \pi'$ iff each block of $\pi$ is contained in a block of $\pi'$. 
\footnote{
NB. The lattice $\Pi_{\le n}$ has a different partial order than the Dowling lattice $Q_n(1) \cong \Pi_{n+1}$ detailed in \cite{Dowling}: by way of comparison, in $Q_n$, $\pi \le \pi'$ iff each block of $\pi'$ is a union of blocks of $\pi$.
}
For economy of notation, we shall write $|\pi|$ for the number of blocks of $\pi \in \Pi_{\le n}$.
Define $F_{\le n} : \Pi_{\le n} \rightarrow \Pi_{n+1}$ as follows: if $\sigma = \sigma^{(1)} | \dots | \sigma^{(|\sigma|)} \in \Pi_{\le n}$ and $[n+1] \setminus \bigcup_{j=1}^{|\sigma|} \sigma^{(j)} = \{s_k\}_{k \in [r]}$, then $F_{\le n}(\sigma) := \sigma^{(1)} | \dots | \sigma^{(|\sigma|)} | \{s_1\} | \dots | \{s_r\}$. Similarly, define $F_{n+1} : \Pi_{n+1} \rightarrow \Pi_{\le n}$ as follows: if $\tau = \tau^{(1)} | \dots | \tau^{(|\tau|)} \in \Pi_{n+1}$ and $n+1 \in \tau^{(u_+)}$, then $F_{n+1}(\tau) :=  \tau^{(1)} | \dots | \tau^{(u_+)} \setminus \{n+1\} | \dots | \tau^{(|\tau|)}$. An important aspect of the relationship between $\Pi_{\le n}$ and $\Pi_{n+1}$ is captured by the following 

\

\textsc{Proposition.} The pair $(F_{\le n}, F_{n+1})$ is a (monotone) Galois connection. $\Box$

\

Define the \emph{support} $\text{supp } \pi$ of $\pi \in \Pi_{\le n}$ as the union of its blocks. Because $\text{supp } \pi_{(\ell,L)} = \ell^{-1}(L)$, the induced colored digraph $D \langle \text{supp } \pi_{(\ell,L)} \rangle$ is well-defined and $\pi_{(\ell,L)} = \pi_{(\ell,L)}^{(1)} | \dots | \pi_{(\ell,L)}^{(|\pi_{(\ell,L)}|)}$ determines a colored digraph via
\begin{equation}
\label{eq:DigraphInducedByPartialPartition}
D_{\langle \ell,L \rangle} := D \langle \text{supp } \pi_{(\ell,L)} \rangle / \{ \pi_{(\ell,L)}^{(1)}, \dots, \pi_{(\ell,L)}^{(|\pi_{(\ell,L)}|)} \}.
\end{equation}
The intuition behind \eqref{eq:DigraphInducedByPartialPartition} is simply that vertices in the $j$th block $\pi_{(\ell,L)}^{(j)}$ are identified (note that the order in which these identifications take place is immaterial, and that the resulting coloring will generally not be canonical).

\

\textsc{Definition}. Call $D_{\langle \ell,L \rangle}$ the \emph{vertex abstraction} of $D$ with respect to $L$.

\

\textsc{Example}. Consider the following colored digraph $D$: 
\begin{center}
	\begin{tikzpicture}[scale=1.25,->,>=stealth',shorten >=1pt]
		\node [draw,circle,minimum size=5mm] (v01) at (0,0) {$1$};
		\node [draw,circle,minimum size=5mm] (v02) at (1,0) {$2$};
		\node [draw,circle,minimum size=5mm] (v03) at (2,0) {$1$};
		\node [draw,circle,minimum size=5mm] (v04) at (3,0) {$3$};
		\foreach \from/\to in {
			v01/v02, v02/v03, v03/v04}
			\draw (\from) to (\to);
	\end{tikzpicture}
\end{center}
We have the following table:

\begin{center}
	\bgroup
	\def\arraystretch{1.5}%  1 is the default, change whatever you need
	\begin{tabular}{ | c | c | c | c | c | c | c | c | c | }
		\hline 
		$L$ & $\varnothing$ & $\{1\}$ & $\{2\}$ & $\{3\}$ & $\{1,2\}$ & $\{1,3\}$ & $\{2,3\}$ & $\{1,2,3\}$
		\\ \hline
		$D_{\langle \ell,L \rangle}$ &
			\ & 
			\begin{tikzpicture}[scale=1.25,->,>=stealth',shorten >=1pt]
				\node [draw,color=white,circle,minimum size=5mm] (v00) at (0,0.2) {};	% vertical padding hack
				\node [draw,circle,minimum size=5mm] (v01) at (0,0) {$1$};
			\end{tikzpicture} &
			\begin{tikzpicture}[scale=1.25,->,>=stealth',shorten >=1pt]
				\node [draw,circle,minimum size=5mm] (v02) at (1,0) {$2$};
			\end{tikzpicture} &
			\begin{tikzpicture}[scale=1.25,->,>=stealth',shorten >=1pt]
				\node [draw,circle,minimum size=5mm] (v03) at (2,0) {$3$};
			\end{tikzpicture} &
			\begin{tikzpicture}[scale=1.25,->,>=stealth',shorten >=1pt]
				\node [draw,circle,minimum size=5mm] (v01) at (0,0) {$1$};
				\node [draw,circle,minimum size=5mm] (v02) at (1,0) {$2$};
				\foreach \from/\to in {
					v01/v02, v02/v01}
					\draw (\from) to (\to);
			\end{tikzpicture} &
			\begin{tikzpicture}[scale=1.25,->,>=stealth',shorten >=1pt]
				\node [draw,circle,minimum size=5mm] (v01) at (0,0) {$1$};
				\node [draw,circle,minimum size=5mm] (v04) at (1,0) {$3$};
				\foreach \from/\to in {
					v01/v04}
					\draw (\from) to (\to);
			\end{tikzpicture} &
			\begin{tikzpicture}[scale=1.25,->,>=stealth',shorten >=1pt]
				\node [draw,circle,minimum size=5mm] (v02) at (0,0) {$2$};
				\node [draw,circle,minimum size=5mm] (v04) at (1,0) {$3$};
			\end{tikzpicture} &
			\begin{tikzpicture}[scale=1.25,->,>=stealth',shorten >=1pt]
				\node [draw,circle,minimum size=5mm] (v02) at (0,0) {$2$};
				\node [draw,circle,minimum size=5mm] (v01) at (1,0) {$1$};
				\node [draw,circle,minimum size=5mm] (v03) at (2,0) {$3$};
				\foreach \from/\to in {
					v01/v02, v02/v01, v01/v03}
					\draw (\from) to (\to);
			\end{tikzpicture}
		\\ \hline
	\end{tabular} $\Box$
	\egroup
\end{center} 

It is easy to see that the map $\pi_{(\ell,\cdot)} : 2^{\Lambda} \rightarrow \Pi_{\le n}$ on the subset lattice $2^{\Lambda}$ is monotone, i.e., if $L \subseteq L'$, then $\pi_{(\ell,L)} \le \pi_{(\ell,L')}$. Furthermore, there is a surjective morphism $\chi_{\ell,L',L} : D_{\langle \ell,L \rangle} \rightarrow D_{\langle \ell,L' \rangle}$ of colored digraphs defined by contracting the preimages of each element of $L' \setminus L$, and $\chi_{\ell,L'',L'} \circ \chi_{\ell,L',L} = \chi_{\ell,L'',L}$. This and some definition-checking yields the following

\

\textsc{Lemma.} When endowed with the obvious refinement morphisms and the $\chi_{\ell,\cdot,\cdot}$, respectively, $\{\pi_{(\ell,L)}\}_{L \subseteq \Lambda}$ and $\{D_{\langle \ell,L \rangle} \}_{L \subseteq \Lambda}$ are both categories. Furthermore, $\pi_{(\ell,\cdot)} : 2^{\Lambda} \rightarrow \{\pi_{(\ell,L)}\}_{L \subseteq \Lambda}$ and $D : \{\pi_{(\ell,L)}\}_{L \subseteq \Lambda} \rightarrow \{D_{\langle \ell,L \rangle} \}_{L \subseteq \Lambda}$ are both functors that yield equivalences of categories. $\Box$

\

In particular, the lattice structure of $2^\Lambda$ is duplicated in $\{\pi_{(\ell,L)}\}_{L \subseteq \Lambda}$ and (since we are considering colored digraphs) $\{D_{\langle \ell,L \rangle} \}_{L \subseteq \Lambda}$. Another noteworthy consequence of this lemma is that the pullback and pushout squares for set inclusions have analogues for partial partitions and vertex abstractions. For example,

\

\textsc{Proposition.} Let $L_1, L_2 \subseteq L$. The pullback of $\chi_{\ell,L,L_1}$ and $\chi_{\ell,L,L_2}$ is given by $\chi_{\ell,L_1,L_1 \cap L_2}$ and $\chi_{\ell,L_2,L_1 \cap L_2}$. $\Box$

\section{\label{sec:DetoursBypassesAbstractions}Detours, bypasses, and path abstractions}

\begin{figure}
	\begin{tikzpicture}[->,>=stealth',shorten >=1pt]
		\node [draw,circle,minimum size=5mm] (v1) at (-2,1) {$1$};
		\node [draw,circle,minimum size=5mm] (v2) at (-2,2) {$2$};
		\node [draw,circle,minimum size=5mm] (v3) at (-1,4) {$3$};
		\node [draw,circle,minimum size=5mm] (v4) at (0,3) {$4$};
		\node [draw,circle,minimum size=5mm] (v5) at (1,4) {$5$};
		\node [draw,circle,minimum size=5mm] (v6) at (2,2) {$6$};
		\node [draw,circle,minimum size=5mm] (v7) at (2,1) {$7$};
		\foreach \from/\to in {
			v1/v4, v2/v4, v4/v6, v4/v7}
		\draw (\from) to (\to);
		\draw (v3) [out=-15,in=105,looseness=1] to (v4);
		\draw (v5) [out=-105,in=15,looseness=1] to (v4);
		\draw (v4) [out=165,in=-75,looseness=1] to (v3);
		\draw (v4) [out=75,in=-165,looseness=1] to (v5);
	\end{tikzpicture}
	\quad \quad \quad
	\begin{tikzpicture}[->,>=stealth',shorten >=1pt]
		\node [draw,circle,minimum size=5mm] (v1) at (-2,1) {$1$};
		\node [draw,circle,minimum size=5mm] (v2) at (-2,2) {$2$};
		\node [draw,circle,minimum size=5mm] (v3) at (-1,4) {$3$};
		\node [draw,circle,minimum size=5mm] (v5) at (1,4) {$5$};
		\node [draw,circle,minimum size=5mm] (v6) at (2,2) {$6$};
		\node [draw,circle,minimum size=5mm] (v7) at (2,1) {$7$};
		\foreach \from/\to in {
			v1/v3, v1/v5, v1/v6, v1/v7, v2/v3, v2/v5, v2/v6, v2/v7, v3/v6, v3/v7, v5/v6, v5/v7}
		\draw (\from) to (\to);
		\draw (v3) [out=30,in=150,looseness=1] to (v5);
		\draw (v5) [out=-150,in=-30,looseness=1] to (v3);
	\end{tikzpicture}
	\caption{\label{fig:LocalStructure} (L) With $D$ the 7-vertex graph shown, we have $V(D)_4^- = \{1,2\}$, $V(D)_4^\pm = \{3,5\}$, and $V(D)_4^+ = \{6,7\}$. (R) $D \upuparrows 4$.} 
\end{figure}
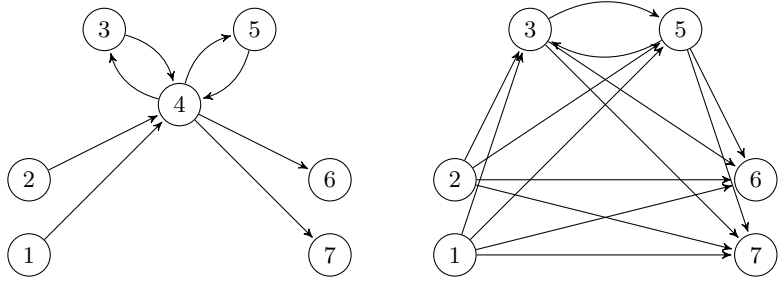

For $v \in V(D)$, define
\begin{eqnarray}
V(D)_v^- & := & \{ x \in V(D) \setminus \{v\} : \mu_D(x,v) \ne 0 \wedge \mu_D(v,x) = 0 \}; \nonumber \\
V(D)_v^\pm & := & \{ x \in V(D) \setminus \{v\} : \mu_D(x,v) \ne 0 \wedge \mu_D(v,x) \ne 0 \}; \nonumber \\
V(D)_v^+ & := & \{ x \in V(D) \setminus \{v\} : \mu_D(x,v) = 0 \wedge \mu_D(v,x) \ne 0 \}; \nonumber \\
V(D)_v^0 & := & \{ x \in V(D) \setminus \{v\} : \mu_D(x,v) = 0 \wedge \mu_D(v,x) = 0 \},
\end{eqnarray}
noting that $\{ \{v\}, V(D)_v^-, V(D)_v^\pm, V(D)_v^+, V(D)_v^0 \}$ forms a partition of $V(D)$. In particular, $P(D)_v := V(D)_v^- \cup V(D)_v^\pm$ is the set of predecessors of $v$, $S(D)_v := V(D)_v^\pm \cup V(D)_v^+$ is the set of successors of $v$, and $P(D)_v \cap S(D)_v = V(D)_v^\pm$. (See figure \ref{fig:LocalStructure}.)

Define the \emph{detour at $v$} by $D \uparrow v$ by $V(D \uparrow v) := V(D)$ and
\begin{equation}
\label{eq:detour}
\mu_{D \uparrow v}(x,y) := 
\begin{cases}
0 & \text{ if } (x,y) \in [P(D)_v \times \{v\}] \cup [\{v\} \times S(D)_v]; \\
1 & \text{ if } (x,y) \in [P(D)_v \times S(D)_v] \setminus \Delta(V(D)); \\
\mu_D(x,y) & \text{ otherwise.}
\end{cases}
\end{equation}
(Here as usual $\Delta$ denotes the diagonal functor.) That is, the detour at $v$ is formed by deleting all arcs involving $v$ and inserting arcs from every predecessor of $v$ to every distinct successor of $v$. Finally, define the \emph{bypass at $v$} by $D \upuparrows v := (D \uparrow v) - v$. (See figure \ref{fig:LocalStructure}.) By construction we have the following

\

\textsc{Proposition.} If $u$, $v$, and $w$ are distinct elements of $V(D)$ and there is a path in $D$ from $u$ to $w$, then there are also paths in both $D \uparrow v$ and $D \upuparrows v$ from $u$ to $w$. $\Box$

\

\textsc{Lemma.} If $D \upuparrows v$ is strongly connected and $P(D)_v, S(D)_v \ne \varnothing$, then $D$ is strongly connected. 

\

\textsc{Proof.} Let $x,y \in [n] \setminus \{v\}$ and let $\gamma(x,y)$ be a path in $D \upuparrows v$ from $x$ to $y$. If $\gamma(x,y)$ does not contain an arc of the form $(u,w)$ for some $u \in P(D)_v$ and $w \in S(D)_v$, then it lifts to a path from $x$ to $y$ in $D$, so assume otherwise. Now $\gamma(x,y)$ is a concatenation of paths of the form $\gamma(x,u) \gamma(u,w) \gamma(w,y)$, which corresponds to a path concatenation of the form $\gamma(x,u) \gamma(u,v) \gamma(v,w) \gamma(w,y)$ in $D$. $\Box$

\

It is easy to see that $D \uparrow v$ and $D \upuparrows v$ both contain a complete digraph with vertex set $V(D)_v^\pm$: thus if this set has more than a single element, the detour and bypass are necessarily cyclic, and their transitive reductions may not be unique. Similarly, there are cases where $D \uparrow v$ has multiple Hamilton cycles (which are also transitive reductions). Therefore in general there is no unique \emph{minimal} digraph with the path preservation property of the preceding proposition. On the other hand, we have the following

\

\textsc{Proposition.} If $D$ is acyclic (so that in particular $V(D)_v^\pm \equiv \varnothing$), then so are $D \uparrow v$ and $D \upuparrows v$. Furthermore, in this event the transitive reductions of $D \uparrow v$ and $D \upuparrows v$ are the unique minimal digraphs on the respective vertex sets $V(D)$ and $V(D) \setminus \{v\}$ such that if $u$, $v$, and $w$ are distinct elements of $V(D)$ and there is a path in $D$ from $u$ to $w$, then there are also paths in both $D \uparrow v$ and $D \upuparrows v$ from $u$ to $w$. $\Box$

\

There are cases where $D \uparrow v$ and $D \upuparrows v$ are their own transitive reductions: e.g., consider a digraph $D$ with only the three arcs $(u,v)$, $(v,w)$, and $(u,w)$: the only arc in $D \uparrow v$ is $(u,w)$. With this in mind, there is a sense in which the detour and bypass can be considered optimal (though typically not minimal) constructions with respect to path preservation in generic digraphs.

\

\textsc{Lemma.} $(D \uparrow v) \uparrow w = (D \uparrow w) \uparrow v$.

\

\textsc{Proof.} See \S \ref{sec:CommutativeDetours} for a naive case analysis. We will prove a more general result in \S \ref{sec:WeightedDigraphs} in a much more elegant and concise manner. $\Box$

\

\textsc{Theorem.} For $U = \{u_j\}_{j \in [m]} \subseteq V(D)$ the obvious generalizations of detour $D \uparrow U := ( \dots (D \uparrow u_1) \dots \uparrow u_{m-1}) \uparrow u_m$ and bypass $D \upuparrows U := (D \uparrow U) - U$ are well-defined. $\Box$

\

Surprisingly, the only reference we could find that even suggests the detour/bypass constructions is \cite{NaumovichClarkeCobleigh}, which seems to take the preceding theorem for granted.

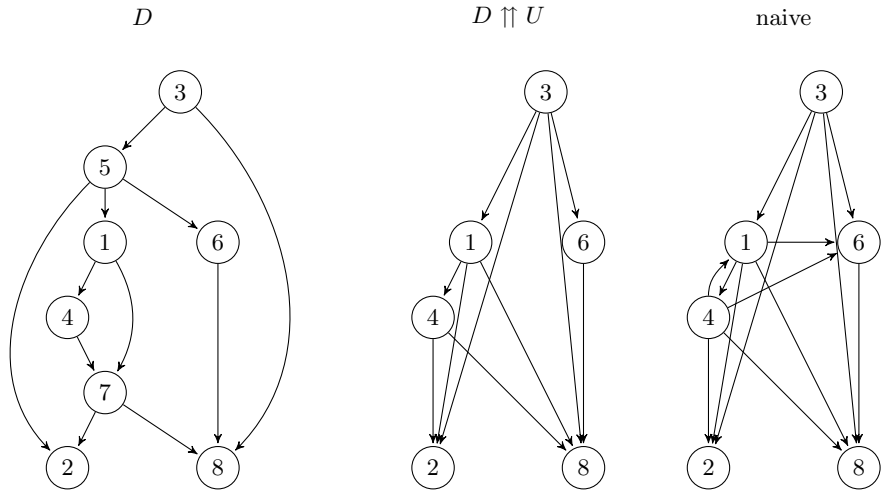
\begin{figure}
	\begin{tikzpicture}[scale=1,->,>=stealth',shorten >=1pt]
		\node at (3,7) {$D$};
		\node [draw,circle,minimum size=5mm] (v1) at (2.5,4) {$1$};
		\node [draw,circle,minimum size=5mm] (v2) at (2,1) {$2$};
		\node [draw,circle,minimum size=5mm] (v3) at (3.5,6) {$3$};
		\node [draw,circle,minimum size=5mm] (v4) at (2,3) {$4$};
		\node [draw,circle,minimum size=5mm] (v5) at (2.5,5) {$5$};
		\node [draw,circle,minimum size=5mm] (v6) at (4,4) {$6$};
		\node [draw,circle,minimum size=5mm] (v7) at (2.5,2) {$7$};
		\node [draw,circle,minimum size=5mm] (v8) at (4,1) {$8$};
		\foreach \from/\to in {
			v1/v4, v3/v5, v4/v7, v5/v1, v5/v6, v6/v8, v7/v2, v7/v8}
		\draw (\from) to (\to);
		\draw (v1) [out=-60,in=60,looseness=1] to (v7);
		\draw (v3) [out=-45,in=45,looseness=1] to (v8);
		\draw (v5) [out=-135,in=135,looseness=1] to (v2);
	\end{tikzpicture}
	\quad \quad \quad
	\begin{tikzpicture}[scale=1,->,>=stealth',shorten >=1pt]
		\node at (3,7) {$D \upuparrows U$};
		\node [draw,circle,minimum size=5mm] (v1) at (2.5,4) {$1$};
		\node [draw,circle,minimum size=5mm] (v2) at (2,1) {$2$};
		\node [draw,circle,minimum size=5mm] (v3) at (3.5,6) {$3$};
		\node [draw,circle,minimum size=5mm] (v4) at (2,3) {$4$};
		\node [draw,circle,minimum size=5mm] (v6) at (4,4) {$6$};
		\node [draw,circle,minimum size=5mm] (v8) at (4,1) {$8$};
		\foreach \from/\to in {
			v1/v2, v1/v4, v1/v8, v3/v1, v3/v2, v3/v6, v3/v8, v4/v2, v4/v8, v6/v8}
		\draw (\from) to (\to);
	\end{tikzpicture}
	\quad \quad \quad
	\begin{tikzpicture}[scale=1,->,>=stealth',shorten >=1pt]
		\node at (3,7) {naive};
		\node [draw,circle,minimum size=5mm] (v1) at (2.5,4) {$1$};
		\node [draw,circle,minimum size=5mm] (v2) at (2,1) {$2$};
		\node [draw,circle,minimum size=5mm] (v3) at (3.5,6) {$3$};
		\node [draw,circle,minimum size=5mm] (v4) at (2,3) {$4$};
		\node [draw,circle,minimum size=5mm] (v6) at (4,4) {$6$};
		\node [draw,circle,minimum size=5mm] (v8) at (4,1) {$8$};
		\foreach \from/\to in {
			v1/v2, v1/v4, v1/v6, v1/v8, v3/v1, v3/v2, v3/v6, v3/v8, v4/v2, v4/v6, v4/v8, v6/v8}
		\draw (\from) to (\to);
		\draw (v4) [out=90,in=-135,looseness=1] to (v1);
	\end{tikzpicture}
	\caption{\label{fig:bypass} (Left) A digraph $D$. (Center) $D \upuparrows U$ for $U = \{5,7\}$. (Right) The result of removing all arcs involving a vertex in $U$, then inserting arcs from every external predecessor of a vertex in $U$ to every distinct external successor of a vertex in $U$ before removing $U$ itself. Note the spurious arcs $(1,6)$, $(4,1)$ and $(4,6)$ that result from this naive procedure.}
\end{figure}

Note that the construction of $D \uparrow U$ is not so simple as removing all arcs involving a vertex in $U$, then inserting arcs from every external predecessor of a vertex in $U$ to every distinct external successor of a vertex in $U$. For example, consider $D$ given by the path of length 3, i.e. $D =
	\begin{tikzpicture}[scale=.5,->,>=stealth',shorten >=1pt]
		\coordinate (v1) at (1,0);
		\coordinate (v2) at (2,0);
		\coordinate (v3) at (3,0);
		\coordinate (v4) at (4,0);
		\draw (v1) [out=30,in=150,looseness=1] to (v2);
		\draw (v2) [out=30,in=150,looseness=1] to (v3);
		\draw (v3) [out=30,in=150,looseness=1] to (v4);
	\end{tikzpicture}$
and $U$ the set whose members are the source and target vertices of $D$. Then $D \upuparrows U = 
	\begin{tikzpicture}[scale=.5,->,>=stealth',shorten >=1pt]
		\coordinate (v2) at (2,0);
		\coordinate (v3) at (3,0);
		\draw (v2) [out=30,in=150,looseness=1] to (v3);
	\end{tikzpicture}$
while the naive procedure mentioned just above yields $
	\begin{tikzpicture}[scale=.5,->,>=stealth',shorten >=1pt]
		\coordinate (v2) at (2,0);
		\coordinate (v3) at (3,0);
		\draw (v2) [out=30,in=150,looseness=1] to (v3);
		\draw (v3) [out=-150,in=-30,looseness=1] to (v2);
	\end{tikzpicture}$. Another example is shown in figure \ref{fig:bypass}.

\

\textsc{Corollary.} If $U \subseteq V(D)$ and $D$ is acyclic, then so are $D \uparrow U$ and $D \upuparrows U$. Furthermore, in this event the transitive reductions of $D \uparrow U$ and $D \upuparrows U$ are the unique minimal digraphs on the respective vertex sets $V(D)$ and $V(D) \setminus U$ such that if $v,w \in V(D) \setminus U$ are distinct and there is a path in $D$ from $v$ to $w$, then there are also paths in both $D \uparrow U$ and $D \upuparrows U$ from $v$ to $w$. $\Box$

\

Recall that for a directed pseudograph or quiver $Q$, the \emph{free category} $F(Q)$ is the category with objects given by vertices of $Q$ and morphisms given by paths in $Q$, with composition given by path concatenation.

\

\textsc{Proposition.} For $U \subseteq V(D)$, there is a functor $F(D) \rightarrow F(D \uparrow U)$ defined on objects as the identity map and on morphisms as the map which deletes elements of $U$ from paths. $\Box$
\footnote{
As suggested in \S \ref{sec:Introduction}, one might consider using this proposition as the foundation of the $\upuparrows$ operator: i.e., given a language $\mathcal{L}$ over $[n]$ whose elements correspond to paths in some digraph $D$, define $\mathcal{L} \upuparrows U$ to be the language over $[n] \setminus U$ obtained by deleting all symbols in $U$ from elements of $\mathcal{L}$. However, it is not obvious from this essentially automata-theoretical persepective that the elements of $\mathcal{L} \upuparrows U$ themselves correspond to paths in some (here) notional digraph $D \upuparrows U$. It seems that a path coherence result of the type needed to establish such a fact would be tantamount to the preceding theorem.
}

\

\textsc{Example.} The digraphs $D$ and $D \upuparrows \{5,7\}$ depicted in figure \ref{fig:bypass} are DAGs. A quick calculation shows that each has 7 paths from a source to a target: the correspondence between them is shown in the table below.
\begin{center}
	\bgroup
	\begin{tabular}{ | c | c | }
		\hline 
		path in $D$ & corresponding path in $D \upuparrows \{5,7\}$
		\\ \hline
		$3 \rightarrow 5 \rightarrow 2$ & $3 \rightarrow 2$ \\
		$3 \rightarrow 5 \rightarrow 1 \rightarrow 7 \rightarrow 2$ & $3 \rightarrow 1 \rightarrow 2$ \\
		$3 \rightarrow 5 \rightarrow 1 \rightarrow 4 \rightarrow 7 \rightarrow 2$ & $3 \rightarrow 1 \rightarrow 4 \rightarrow 2$ \\
		$3 \rightarrow 8$ & $3 \rightarrow 8$ \\
		$3 \rightarrow 5 \rightarrow 6 \rightarrow 8$ & $3 \rightarrow 6 \rightarrow 8$ \\
		$3 \rightarrow 5 \rightarrow 1 \rightarrow 7 \rightarrow 8$ & $3 \rightarrow 1 \rightarrow 8$ \\
		$3 \rightarrow 5 \rightarrow 1 \rightarrow 4 \rightarrow 7 \rightarrow 8$ & $3 \rightarrow 1 \rightarrow 4 \rightarrow 8$ \\
		\hline
	\end{tabular} $\Box$
	\egroup
\end{center} 

\

\begin{figure}
	\begin{tikzpicture}[->,>=stealth',shorten >=1pt]
		\node [draw,circle,minimum size=7mm] (v1) at (0,0) {1};
		\node [draw,circle,minimum size=7mm] (v2) at (2,0) {2};
		\node [draw,circle,minimum size=7mm] (v3) at (2,2) {3};
		\node [draw,circle,minimum size=7mm] (v4) at (0,2) {4};
		\draw (v1) [out=-30,in=-150,looseness=1] to (v2);
		\draw (v1) [out=30,in=-120,looseness=1] to (v3);
		\draw (v2) [out=60,in=-60,looseness=1] to (v3);
		\draw (v3) [out=-150,in=60,looseness=1] to (v1);
		\draw (v3) [out=150,in=30,looseness=1] to (v4);
		\draw (v4) [out=-120,in=120,looseness=1] to (v1);
	\end{tikzpicture}
	\quad \quad \quad
	\begin{tikzpicture}[->,>=stealth',shorten >=1pt]
		\node [draw,circle,minimum size=7mm] (v1) at (0,0) {1};
		\node [draw,circle,minimum size=7mm] (v2) at (2,0) {2};
		\node [draw,circle,minimum size=7mm] (v4) at (0,2) {4};
		\draw (v1) [out=-30,in=-150,looseness=1] to (v2);
		\draw (v1) [out=60,in=-60,looseness=1] to (v4);
		\draw (v2) [out=150,in=30,looseness=1] to (v1);
		\draw (v2) [out=120,in=-30,looseness=1] to (v4);
		\draw (v4) [out=-120,in=120,looseness=1] to (v1);
	\end{tikzpicture}
	\quad \quad \quad
	\begin{tikzpicture}[->,>=stealth',shorten >=1pt]
		\node [draw,circle,minimum size=7mm] (v1) at (0,0) {1};
		\node [draw,circle,minimum size=7mm] (v2) at (2,0) {2};
		\draw (v1) [out=-30,in=-150,looseness=1] to (v2);
		\draw (v2) [out=150,in=30,looseness=1] to (v1);
	\end{tikzpicture}
	\caption{\label{fig:CyclicExample} (Left) $D$. (Center) $D \upuparrows 3$. (Right) $D \upuparrows 1|24$ (see \eqref{eq:PathAbstraction} for this construction).}
\end{figure}
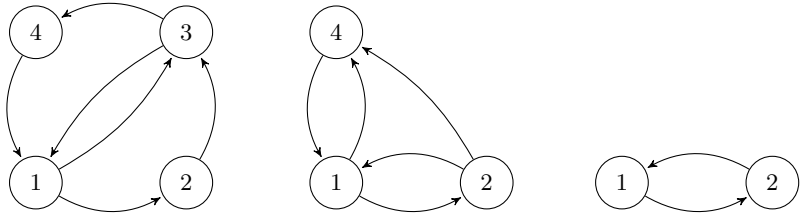

\textsc{Example.} Consider the digraph $D$ in figure \ref{fig:CyclicExample}. There are four interesting cycles: $1 \rightarrow 3 \rightarrow 1 \equiv 131$ (omitting arrows for concision), $1231$, $1341$, and $12341$. These are respectively mapped under $\upuparrows 3$ to $1$ (not a cycle!), $121$, $141$, and $1241$. Subsequently contracting vertices $2$ and $4$ \emph{\`a la} \eqref{eq:PathAbstraction} maps the remaining cycles in turn to the single cycle $121$. This extends to a mapping on all paths, e.g. the path $1231341234$ maps under $\upuparrows 3$ to $1214124$ and subsequently under the contraction of vertices $2$ and $4$ to $121212$. $\Box$

\

From the proposition we see that for $L \subseteq L'$ there is a surjective morphism $\phi_{\ell,L,L'} : F(D \uparrow \ell^{-1}(L)) \rightarrow F(D \uparrow \ell^{-1}(L'))$ defined by sending vertices of $D$ to themselves and deleting elements of $L' \setminus L$ from paths. Therefore we get the following 

\

\textsc{Lemma.} When endowed with the morphisms $\phi_{\ell,\cdot,\cdot}$, $\{F(D \uparrow \ell^{-1}(L))\}_{L \subseteq \Lambda}$ is a category and the map $F_{D,\ell} : L \mapsto F(D \uparrow \ell^{-1}(L))$ is a functor yielding an equivalence of categories. $\Box$

\

We close this section by showing that detour/bypass and contraction operations on disjoint vertex sets commute.

\

\textsc{Lemma.} If $u, v, w \in V(D)$ are distinct, then $(D \uparrow u) / \{v,w\} = (D/\{v,w\}) \uparrow u$.

\

\textsc{Proof.} See \S \ref{sec:CommutativeContractions} for a naive case analysis. We will prove a more general result in \S \ref{sec:WeightedDigraphs} in a much more elegant and concise manner. $\Box$

\

\textsc{Theorem.} The detour and bypass operations commute with disjoint graph contractions: if $U, V_1, \dots, V_m \subseteq V(D)$ are disjoint, then $(D \uparrow U) / \{V_1,\dots,V_m\} = (D / \{V_1,\dots,V_m\}) \uparrow U$, and similarly for bypasses. $\Box$

\

Thus we can freely make the following

\

\textsc{Definition.} The \emph{path abstraction} of $D$ with respect to $\pi \in \Pi_{\le |V(D)|}$ is 
\begin{equation}
\label{eq:PathAbstraction}
D \upuparrows \pi := (D \upuparrows (V(D) \setminus \text{supp } \pi)) / \{\pi^{(1)},\dots,\pi^{(|\pi|)}\} = (D / \{\pi^{(1)},\dots,\pi^{(|\pi|)}\}) \upuparrows (V(D) \setminus \text{supp } \pi).
\end{equation}
We shall write $D \upuparrows (\ell,L) := D \upuparrows \pi_{(\ell,L)}$.

\

Unfortunately, for $L \subseteq L'$ there is not a simple relationship between $D \upuparrows (\ell, L)$ and $D \upuparrows (\ell, L')$. To see this, let $L \subset \Lambda$, let $\ell_0 \in \Lambda \setminus L$, and let $L' := L \cup \{\ell_0\}$. If $\pi_{(\ell,L)} = \pi_{(\ell,L)}^{(1)} | \dots | \pi_{(\ell,L)}^{(|\pi_{(\ell,L)}|)}$ and $\pi_{(\ell,L')} = \pi_{(\ell,L')}^{(1)} | \dots | \pi_{(\ell,L')}^{(|\pi_{(\ell,L')}|)} = \pi_{(\ell,L)}^{(1)} | \dots | \pi_{(\ell,L)}^{(|\pi_{(\ell,L)}|)} | \pi_{(\ell,L')}^{(|\pi_{(\ell,L)}|+1)}$, defining $ \pi_{(\ell,L')}^{(0)} := V(D) \setminus \bigcup_{j \in [|\pi_{(\ell,L)}|+1]} \pi_{(\ell,L')}^{(j)}$ yields a partition $\{\pi_{(\ell,L')}^{(j)}\}_{j = 0}^{|\pi_{(\ell,L)}|+1}$ of $V(D)$, and it is readily seen that $D \upuparrows (\ell, L) = ((D \upuparrows \pi_{(\ell,L)}^{(0)}) / \{\pi_{(\ell,L)}^{(1)},\dots,\pi_{(\ell,L)}^{(|\pi_{(\ell,L)}|)}\}) \upuparrows \pi_{(\ell,L')}^{(|\pi_{(\ell,L)}|+1)}$ while $D \upuparrows (\ell,L') = ((D \upuparrows \pi_{(\ell,L)}^{(0)}) / \{\pi_{(\ell,L)}^{(1)},\dots,\pi_{(\ell,L)}^{(|\pi_{(\ell,L)}|)}\}) \setminus \pi_{(\ell,L')}^{(|\pi_{(\ell,L)}|+1)}$. That is, the essential distinction between $D \upuparrows (\ell,L)$ and $D \upuparrows (\ell,L')$ is that the former has a bypass while the latter has a contraction. These two operations are typically not readily comparable, and so we defer the quest for a structure theory of path abstractions.

\section{\label{sec:WeightedDigraphs}Weighted digraphs}

Generalizing the constructions of \S \ref{sec:DetoursBypassesAbstractions} to weighted digraphs introduces some subtleties. However, it also leads to simpler proofs.

As a preliminary step, consider the case of directed multigraphs. For a directed multigraph $D$, there is a unique minimal subdivision $D^\oslash$ of $D$ that is a digraph. Thus for $v \in V(D)$, $D^\oslash \uparrow v$ is well-defined. For $x,y \in V(D)$, write $\nu_D(x,y)$ for the number of walks in $D$ from $x$ to $y$, so that $\nu_D(x,y) < \infty$ iff $x$ and $y$ are not in the same strong component of $D$ (of course, this is automatically the case if $D$ is acyclic [which also implies that $D^\oslash$ and $D^\oslash \uparrow v$ are acyclic]). For $x,y \in V(D)$, it is clear that $\nu_{D}(x,y) = \nu_{D^\oslash}(x,y)$, and moreover that any reasonable definition of $D \uparrow v$ should satisfy $\nu_{D^\oslash \uparrow v}(x,y) = \nu_{D \uparrow v}(x,y)$. Meanwhile, note that \eqref{eq:detour} is equivalent to
\begin{equation}
\mu_{D \uparrow v}(x,y) = 
\begin{cases}
0 & \text{ if } (x,y) \in [V(D) \times \{v\}] \cup [\{v\} \times V(D)] \cup \Delta(V(D)); \\
\mu_D(x,y) \vee (\mu_D(x,v) \wedge \mu_D(v,y)) & \text{ otherwise.}
\end{cases} \nonumber
\end{equation}

These considerations indicate that for directed multigraphs, we should simply replace $\vee$ with $+$ and $\wedge$ with $\cdot$ \emph{\`a la}
\begin{equation}
\label{eq:DAMGDetour}
\mu_{D \uparrow v}(x,y) :=
\begin{cases}
0 & \text{ if } (x,y) \in [V(D) \times \{v\}] \cup [\{v\} \times V(D)] \cup \Delta(V(D)); \\
\mu_D(x,y) + \mu_D(x,v) \cdot \mu_D(v,y) & \text{ otherwise.}
\end{cases}
\end{equation}
To generalize further from directed multigraphs to weighted digraphs, the addition and multipication operations above can be taken to be those of a commutative semiring that the weights are presumed to inhabit, and $\mu_D$ can be taken to indicate the weights (or for the further generalization of a weighted directed multigraph, the appropriate sum of weights). However, while the RHS of \eqref{eq:DAMGDetour} is always well-defined, in many circumstances it leads to behavior that is more troublesome than for the special case of unweighted digraphs. 

For convenience, we shall write, e.g., $\mu_{xy} := \mu_D(x,y)$ in the rest of this section.

\

\textsc{Theorem.} $(D \uparrow v) \uparrow w \ne (D \uparrow w) \uparrow v$ iff $\mu_{vw} \mu_{wv} \ne 0$ and there exists $(x,y) \in V(D)^2 \setminus ([V(D) \times \{v,w\}] \cup [\{v,w\} \times V(D)] \cup \Delta(V(D)))$ such that $\mu_{xv} \mu_{vy} \ne \mu_{xw} \mu_{wy}$. 

\

\textsc{Proof.} We have that
\begin{equation}
\label{eq:DAMGDoubleDetour}
\mu_{(D \uparrow v) \uparrow w}(x,y) :=
\begin{cases}
0 & \text{ if } (x,y) \in [V(D) \times \{v,w\}] \cup [\{v,w\} \times V(D)] \cup \Delta(V(D)); \\
\mu_{D \uparrow v}(x,y) + \mu_{D \uparrow v}(x,w) \cdot \mu_{D \uparrow v}(w,y) & \text{ otherwise.}
\end{cases} 
\end{equation}
Applying \eqref{eq:DAMGDetour} twice shows that the nonzero entries of $\mu_{(D \uparrow v) \uparrow w}(x,y)$ are of the form
\begin{equation}
\mu_{xy} + \mu_{xv} \mu_{vy} + \mu_{xw} \mu_{wy} + \mu_{xv} \mu_{vw} \mu_{wy} + \mu_{xw} \mu_{wv} \mu_{vy} + \mu_{xv} \mu_{vw} \mu_{wv} \mu_{vy}. \nonumber
\end{equation}
The expression above is not symmetric under the exchange of $v$ and $w$ owing purely to the last term (note that in the special case of digraphs addressed in \S \ref{sec:DetoursBypassesAbstractions}, the concomitant Boolean semiring operations recover this lost symmetry as required since the last term is nonzero only if the second term is also). The conditions under which there exist $x,y$ such that $x \ne y$ and $\mu_{xv} \mu_{vw} \mu_{wv} \mu_{vy} \ne \mu_{xw} \mu_{wv} \mu_{vw} \mu_{wy}$ are stated in the theorem. $\Box$

\

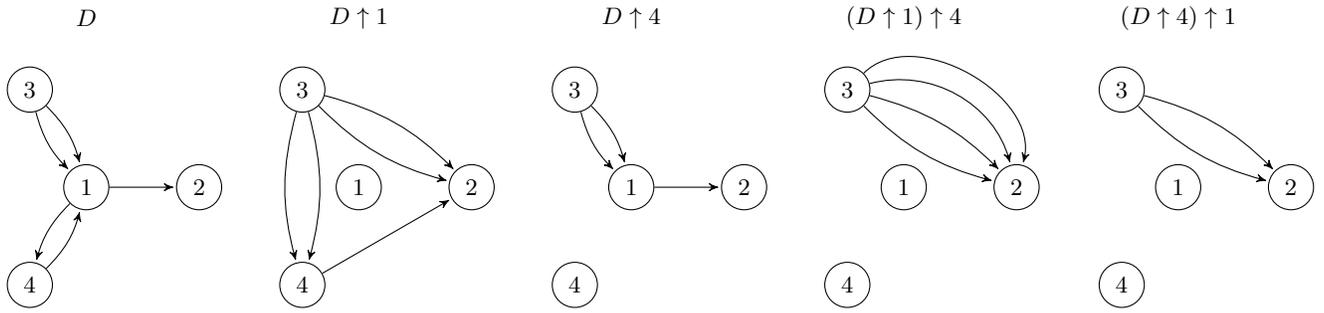
\begin{figure}
	\begin{tikzpicture}[scale=0.75,->,>=stealth',shorten >=1pt]
		\node at (0,3) {$D$};
		\node [draw,circle,minimum size=6mm] (v1) at (0,0) {1};
		\node [draw,circle,minimum size=6mm] (v2) at (2,0) {2};
		\node [draw,circle,minimum size=6mm] (v3) at (-1,1.7321) {3};
		\node [draw,circle,minimum size=6mm] (v4) at (-1,-1.7321) {4};
		\draw (v1) to (v2);
		\draw (v1) [out=-135,in=75,looseness=1] to (v4);
		\draw (v3) [out=-45,in=105,looseness=1] to (v1);
		\draw (v3) [out=-75,in=135,looseness=1] to (v1);
		\draw (v4) [out=45,in=-105,looseness=1] to (v1);
	\end{tikzpicture}
	\quad \quad 
	\begin{tikzpicture}[scale=0.75,->,>=stealth',shorten >=1pt]
		\node at (0,3) {$D \uparrow 1$};
		\node [draw,circle,minimum size=6mm] (v1) at (0,0) {1};
		\node [draw,circle,minimum size=6mm] (v2) at (2,0) {2};
		\node [draw,circle,minimum size=6mm] (v3) at (-1,1.7321) {3};
		\node [draw,circle,minimum size=6mm] (v4) at (-1,-1.7321) {4};
		\draw (v3) [out=-15,in=135,looseness=1] to (v2);
		\draw (v3) [out=-45,in=165,looseness=1] to (v2);
		\draw (v3) [out=-75,in=75,looseness=1] to (v4);
		\draw (v3) [out=-105,in=105,looseness=1] to (v4);
		\draw (v4) to (v2);
	\end{tikzpicture}
	\quad \quad 
	\begin{tikzpicture}[scale=0.75,->,>=stealth',shorten >=1pt]
		\node at (0,3) {$D \uparrow 4$};
		\node [draw,circle,minimum size=6mm] (v1) at (0,0) {1};
		\node [draw,circle,minimum size=6mm] (v2) at (2,0) {2};
		\node [draw,circle,minimum size=6mm] (v3) at (-1,1.7321) {3};
		\node [draw,circle,minimum size=6mm] (v4) at (-1,-1.7321) {4};
		\draw (v1) to (v2);
		\draw (v3) [out=-45,in=105,looseness=1] to (v1);
		\draw (v3) [out=-75,in=135,looseness=1] to (v1);
	\end{tikzpicture}
	\quad \quad 
	\begin{tikzpicture}[scale=0.75,->,>=stealth',shorten >=1pt]
		\node at (0,3) {$(D \uparrow 1) \uparrow 4$};
		\node [draw,circle,minimum size=6mm] (v1) at (0,0) {1};
		\node [draw,circle,minimum size=6mm] (v2) at (2,0) {2};
		\node [draw,circle,minimum size=6mm] (v3) at (-1,1.7321) {3};
		\node [draw,circle,minimum size=6mm] (v4) at (-1,-1.7321) {4};
		\draw (v3) [out=45,in=75,looseness=1] to (v2);
		\draw (v3) [out=15,in=105,looseness=1] to (v2);
		\draw (v3) [out=-15,in=135,looseness=1] to (v2);
		\draw (v3) [out=-45,in=165,looseness=1] to (v2);
	\end{tikzpicture}
	\quad \quad 
	\begin{tikzpicture}[scale=0.75,->,>=stealth',shorten >=1pt]
		\node at (0,3) {$(D \uparrow 4) \uparrow 1$};
		\node [draw,circle,minimum size=6mm] (v1) at (0,0) {1};
		\node [draw,circle,minimum size=6mm] (v2) at (2,0) {2};
		\node [draw,circle,minimum size=6mm] (v3) at (-1,1.7321) {3};
		\node [draw,circle,minimum size=6mm] (v4) at (-1,-1.7321) {4};
		\draw (v3) [out=-15,in=135,looseness=1] to (v2);
		\draw (v3) [out=-45,in=165,looseness=1] to (v2);
	\end{tikzpicture}
	\caption{\label{fig:DetoursDontCommute} $(D \uparrow v) \uparrow w \ne (D \uparrow w) \uparrow v$ in general when $D$ is a directed multigraph.}
\end{figure}

\textsc{Example.} Consider the multigraph $D$ shown in figure \ref{fig:DetoursDontCommute}. With $v = 1$ and $w = 4$, the terms $ \mu_{xv} \mu_{vw} \mu_{wv} \mu_{vy}$ and $\mu_{xw} \mu_{wv} \mu_{vw} \mu_{wy}$ respectively correspond to the matrices
\begin{equation}
\begin{pmatrix}0 \\ 0 \\ 2 \\ 1
\end{pmatrix} 
\cdot 
\begin{pmatrix} 1 \end{pmatrix} 
\cdot 
\begin{pmatrix} 1 \end{pmatrix} 
\cdot
\begin{pmatrix}0 & 1 & 0 & 1
\end{pmatrix}
=
\begin{pmatrix}0 & 0 & 0 & 0 \\ 0 & 0 & 0 & 0 \\ 0 & 2 & 0 & 2 \\ 0 & 1 & 0 & 1
\end{pmatrix}
\quad \text{and} \quad
\begin{pmatrix}1 \\ 0 \\ 0 \\ 0
\end{pmatrix} 
\cdot 
\begin{pmatrix} 1 \end{pmatrix} 
\cdot 
\begin{pmatrix} 1 \end{pmatrix} 
\cdot
\begin{pmatrix}1 & 0 & 0 & 0
\end{pmatrix}
=
\begin{pmatrix}1 & 0 & 0 & 0 \\ 0 & 0 & 0 & 0 \\ 0 & 0 & 0 & 0 \\ 0 & 0 & 0 & 0
\end{pmatrix}. \nonumber
\end{equation}
These have entries in $V(D)^2 \setminus ([V(D) \times \{v,w\}] \cup [\{v,w\} \times V(D)] \cup \Delta(V(D))) = \{(2,3),(3,2)\}$ that differ. $\Box$

\

\textsc{Corollary.} If $D$ is acyclic and $U \subseteq V(D)$, then $D \uparrow U$ is well-defined and acyclic. $\Box$

\

There may be situations of practical interest in which $D$ is not acyclic, which necessarily complicates any criterion for establishing the existence of a well-defined detour/bypass. We proceed below to establish the most obvious criterion in this vein.

\

\textsc{Lemma.} If $D$ has no 2-cycles at all and no 3-cycles intersecting $v$, then $D \uparrow v$ has no 2-cycles.

\

\textsc{Proof.} Suppose that $\mu_{D \uparrow v}(x,y) \cdot \mu_{D \uparrow v}(y,x) \ne 0$. By hypothesis, $\mu_{xy} \mu_{yx} = 0$, so without loss of generality assume that $\mu_{xy} = 0$. Now 
\begin{equation}
\mu_{D \uparrow v}(x,y) \cdot \mu_{D \uparrow v}(y,x) = \mu_{xv} \mu_{vy} (\mu_{yx} + \mu_{yv} \mu_{vx}) \ne 0. \nonumber
\end{equation}
In particular, $\mu_{xv} \ne 0 \ne \mu_{vy}$. 

If $\mu_{yx} = 0$, then we must have that $\mu_{xv} \mu_{vx} \ne 0 \ne \mu_{vy} \mu_{yv}$, contradicting the assumption that $D$ has no 2-cycles. If on the other hand $\mu_{yx} \ne 0$, then either $\mu_{yv} \mu_{vx} = 0$ or $\mu_{yv} \mu_{vx} \ne 0$. In the first case, $\mu_{xv} \mu_{vy} \mu_{yx} \ne 0$, contradicting the assumption that $D$ has no 3-cycles intersecting $v$. In the second case, we again contradict the assumption that $D$ has no 2-cycles. Thus it must be that $\mu_{D \uparrow v}(x,y) \cdot \mu_{D \uparrow v}(y,x) = 0$, i.e., $D \uparrow v$ has no 2-cycles. $\Box$

\

While it is tempting to spend the effort to recast the preceding lemma as the base case of an induction, the complexity of finding short cycles in digraphs suggests that any theorem actually resulting from such an exercise would be less useful in practice than simply checking online whether or not successive detours commute. For this reason we elect to move on to the following more useful result:

\

\textsc{Theorem.} If $u, v, w \in V(D)$ are distinct, then $(D \uparrow u) / \{v,w\} = (D/\{v,w\}) \uparrow u$.

\

\textsc{Proof.} Write $X := V(D) \setminus \{u,v,w\}$, so that $\mu_D$ takes the form
\begin{equation}
\bordermatrix{\mu_D & u & v & w & X \cr
	u & 0 & \mu_{uv} & \mu_{uw} & \mu_{uX} \cr
	v & \mu_{vu} & 0 & \mu_{vw} & \mu_{vX} \cr
	w & \mu_{wu} & \mu_{wv} & 0 & \mu_{wX} \cr
	X & \mu_{Xu} & \mu_{Xv} & \mu_{Xw} & \mu_{XX} \cr
	}. \nonumber
\end{equation}
For a generic square matrix $M$ and $n$-tuple $z$, let $d(M)_j := M_{jj}$ and $d(z)_{jk} := z_j \delta_{jk}$, i.e., $d$ is the obvious notion of a diagonal map in both cases. Define $M^\boxbslash := M - d(d(M))$. We then have that $\mu_{D \uparrow u}$ and $\mu_{D / \{v,w\}}$ are respectively
\begin{equation}
\bordermatrix{\mu_{D \uparrow u} & u & v & w & X \cr
	u & 0 & 0 & 0 & 0 \cr
	v & 0 & 0 & \mu_{vw} + \mu_{vu} \mu_{uw} & \mu_{vX} + \mu_{vu} \mu_{uX} \cr
	w & 0 & \mu_{wv} + \mu_{wu} \mu_{uv} & 0 & \mu_{wX} + \mu_{wu} \mu_{uX} \cr
	X & 0 & \mu_{Xv} + \mu_{Xu} \mu_{uv} & \mu_{Xw} + \mu_{Xu} \mu_{uw} & \mu_{XX} + (\mu_{Xu}\mu_{uX})^\boxbslash \cr
	}; \nonumber
\end{equation}
\begin{equation}
\bordermatrix{\mu_{D / \{v,w\}} & u & vw & X \cr
	u & 0 & \mu_{uv} + \mu_{uw} & \mu_{uX} \cr
	vw & \mu_{vu} + \mu_{wu} & 0 & \mu_{vX} + \mu_{wX} \cr
	X & \mu_{Xu} & \mu_{Xv} + \mu_{Xw} & \mu_{XX} \cr
	}. \nonumber
\end{equation}

Consequently, $\mu_{(D \uparrow u) \setminus \{v,w\}}$ and $\mu_{(D / \{v,w\}) \uparrow u}$ are respectively
\begin{equation}
\bordermatrix{\mu_{(D \uparrow u) \setminus \{v,w\}} & u & vw & X \cr
	u & 0 & 0 & 0 \cr
	vw & 0 & 0 & \mu_{vX} + \mu_{vu} \mu_{uX} + \mu_{wX} + \mu_{wu} \mu_{uX}  \cr
	X & 0 & \mu_{Xv} + \mu_{Xu} \mu_{uv} + \mu_{Xw} + \mu_{Xu} \mu_{uw} & \mu_{XX} + (\mu_{Xu}\mu_{uX})^\boxbslash \cr
	}; \nonumber
\end{equation}
\begin{equation}
\bordermatrix{\mu_{(D / \{v,w\}) \uparrow u} & u & vw & X \cr
	u & 0 & 0 & 0 \cr
	vw & 0 & 0 & \mu_{vX} + \mu_{wX} + (\mu_{vu} + \mu_{wu}) \mu_{uX} \cr
	X & 0 & \mu_{Xv} + \mu_{Xw} + \mu_{Xu} (\mu_{uv} + \mu_{uw}) & \mu_{XX} + (\mu_{Xu}\mu_{uX})^\boxbslash \cr
	}. \nonumber
\end{equation}
These are obviously equal. $\Box$

\

\textsc{Corollary.} If $D$ is acyclic, then its path abstraction is well-defined. $\Box$

\section{\label{sec:RandomDigraphs}Random digraphs}

For $0 \le p \le 1$, let $D_{n,p}$ denote the random digraph \cite{FriezeKaronski} with $V(D_{n,p}) = [n]$ and independent probabilities $\mathbb{P}(\mu_{D_{n,p}}(x,y) = 1 | x \ne y) \equiv p$. 

\begin{figure}[htbp]
\includegraphics[trim = 15mm 0mm 15mm 0mm, clip, width=180mm,keepaspectratio]{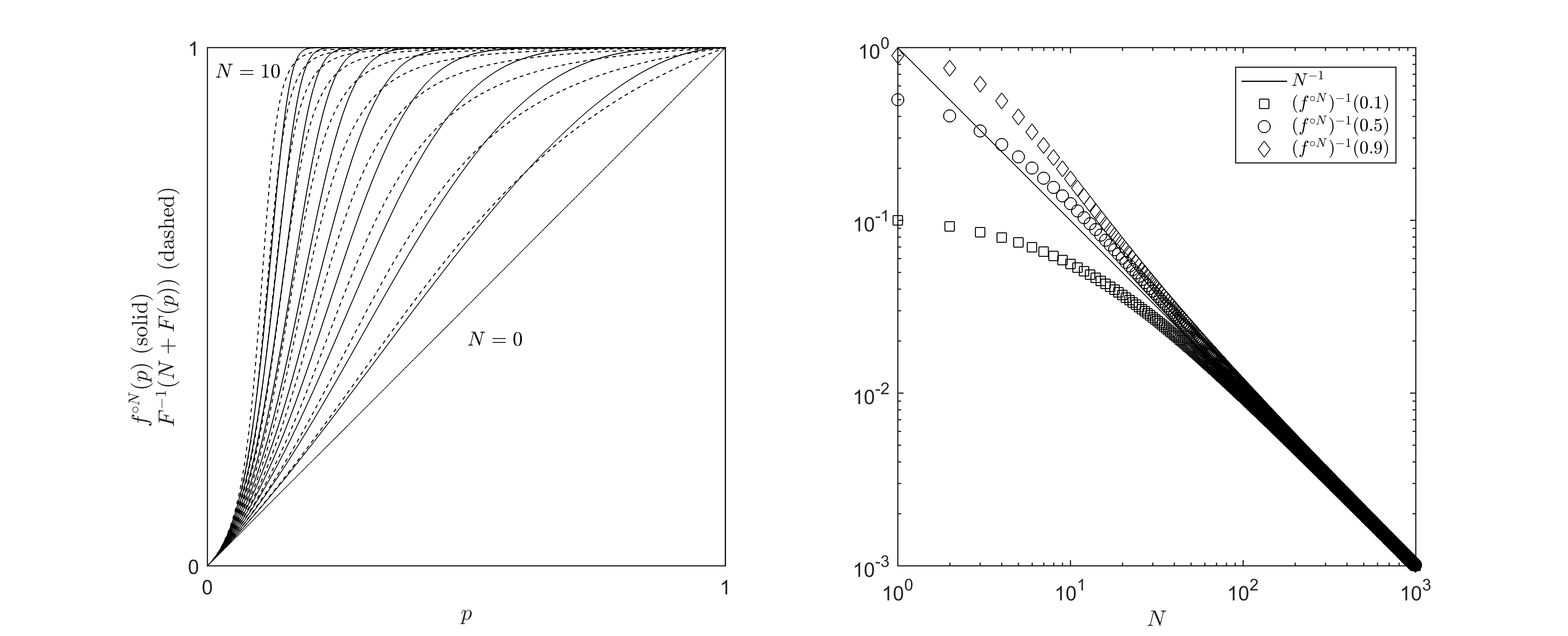}% Here is how to import pix
\caption{ \label{fig:IteratedF2} (L) $f^{\circ N}(p)$ (solid) and $F^{-1}(N+F(p))$ (dashed) for $N \in \{0,\dots,10\}$. (R) Inverse images $(f^{\circ N})^{-1}(p_0)$ for $p_0 \in \{0.1,0.5,0.9\}$.} 
\end{figure}

Let $u \in [n]$. According to \eqref{eq:detour}, there are two ways for the event $\mu_{D_{n,p} \upuparrows u}(x,y) = 1$ to occur for $x \ne y$: either $(x,y) \in P(D_{n,p})_u \times S(D_{n,p})_u$, or $(x,y) \notin P(D_{n,p})_u \times S(D_{n,p})_u$ but $\mu_{D_{n,p}}(x,y) = 1$. These two subevents are disjoint, with the probability of the former equal to $p^2$ and the probability of the latter equal to $(1-p^2)p$, so we have that $\mathbb{P} \left ( \mu_{D_{n,p} \upuparrows u}(x,y) = 1 | x \ne y \right ) = p^2 + (1-p^2)p =: f(p)$. It follows for that for $U \subseteq [n]$
\begin{equation}
\label{eq:RandomBypass1}
\mathbb{P} \left ( \mu_{D_{n,p} \upuparrows U}(x,y) = 1 | x \ne y \right ) = f^{\circ |U|}(p),
\end{equation}
where a $|U|$-fold composition is indicated on the RHS (see figure \ref{fig:IteratedF2}). That is, 
\begin{equation}
\label{eq:RandomBypass2}
D_{n,p} \upuparrows U = D_{n-|U|,f^{\circ |U|}(p)}.
\end{equation}
This suggests the possibility of a renormalization group strategy for studying $D_{n,p}$ (and percolation thresholds in particular), but we limit our discussion of this to a terse remark in \S \ref{sec:Renormalization}.

A qualitative approximation for $f^{\circ N}(p)$ is readily obtained by the following tactic described in \cite{CountIblis}: temporarily writing $p(N) := f^{\circ N}(p)$, we have that $p(N+1)-p(N) = p^2(N) - p^3(N)$. Treating $N$ as a continuous parameter yields the approximation $\frac{dp}{dN} \approx p^2-p^3$. Writing $F(p) := \log \frac{p}{1-p} - \frac{1}{p}$ and noting that $\frac{dF}{dp} = \frac{1}{p^2-p^3}$ yields that 
\begin{equation}
\label{eq:IterateApprox}
f^{\circ N}(p) \approx F^{-1}(N+F(p)).
\end{equation} 
A sophisticated but interesting restatement of \eqref{eq:IterateApprox} is that $F$ is approximately equivariant with respect to the $\mathbb{Z}$-actions on $[0,1]$ and $\mathbb{R}$ given respectively by iterating $f$ and adding.
\footnote{
The approximate $\mathbb{Z}$-equivariance expressed in \eqref{eq:IterateApprox} obviously generalizes to other functions $f$ and $F$ such that $F'(p) = (f(p)-p)^{-1}$. However, we are not aware of results corresponding to this fact or its further generalization to higher-order differences.
}

%It can be shown that the low-order coefficients of $f^{\circ N}(p) =: \sum_{d=1}^\infty a_d(N) \cdot p^d$ are given by 
%\begin{eqnarray}
%a_1(N) & = & 1, \nonumber \\
%a_2(N) & = & N, \nonumber \\
%a_3(N) & = & (N-1)^2-1, \nonumber \\
%a_4(N) & = & (N-1)^3-2(N-1)^2-3(N-1), \nonumber \\
%a_5(N) & = & (N-1)^4-\frac{14}{3}(N-1)^3-2(N-1)^2+\frac{11}{3}(N-1), \nonumber \\
%& \vdots & \nonumber
%\end{eqnarray}
%but we are not aware of a closed form for $a_d(N)$ as a function of $d$, which severely limits the utility of such an expansion. 

Meanwhile, the event $\mu_{D_{n,p}/U}(x,\{U\}) = 1$ occurs for $x \ne \{U\}$ iff there is some $u \in U$ such that $\mu_{D_{n,p}}(x,u) = 1$. Equivalently, the event $\mu_{D_{n,p}/U}(x,\{U\}) = 0$ occurs iff $\mu_{D_{n,p}}(x,u) = 0$ for all $u \in U$, and this event clearly has probability $(1-p)^{|U|}$. Thus 
\begin{equation}
\label{eq:RandomContraction1}
\mathbb{P} \left ( \mu_{D_{n,p}/U}(x,\{U\}) = 1 | x \ne \{U\} \right ) = 1 - (1-p)^{|U|}.
\end{equation}
More generally, if $\pi \in \Pi_n$, then 
\begin{equation}
\label{eq:RandomContraction2}
\mathbb{P} \left ( \mu_{D_{n,p}/\{\pi^{(1)},\dots,\pi^{(|\pi|)}\}}(\{\pi^{(j)}\},\{\pi^{(k)}\}) = 1 | j \ne k \right ) = 1 - (1-p)^{|\pi^{(j)}| \cdot |\pi^{(k)}|}.
\end{equation}

We can combine the preceding observations into the following

\ 

\textsc{Theorem.} If $\pi \in \Pi_{\le n}$, then
\begin{equation}
\label{eq:RandomContraction2}
\mathbb{P} \left ( \mu_{D_{n,p} \upuparrows \pi}(\{\pi^{(j)}\},\{\pi^{(k)}\}) = 1 | j \ne k \right ) = 
1 - \left ( 1-f^{\circ (n-|\text{supp } \pi|)}(p) \right )^{|\pi^{(j)}| \cdot |\pi^{(k)}|}. 
\end{equation}
In particular, the expected number of arcs in $D_{n,p} \upuparrows \pi$ is the sum over $j \ne k$ of the RHS of \eqref{eq:RandomContraction2} (note that the number of vertices is just $|\pi|$). $\Box$

\

	%% NYC street data from http://www.dis.uniroma1.it/challenge9/download.shtml
	%
	%%% Get indices for the heart of the financial district
	%indV = find((lon>=-7.40135e7)&(lon<=-7.40085e7)&...
	%    (lat>=4.07040e7)&(lat<=4.07090e7)&...
	%    (lat-4.07040e7>=(0.47)*(lon+7.40125e7))&...
	%    (lat-4.07065e7<=lon+7.40135e7)&...
	%    (lat-4.07090e7<=(-.60)*(lon+7.40110e7)));
	%%% Produce adjacency matrix
	%% It is easier to produce arcs by hand since the arc orientation data is
	%% bad. We assume (perhaps wrongly) that New St. runs parallel to Broadway
	%% and get the other street directions from Google Maps.
	%D = sparse(28,28);  % numel(indV) = 28
	%A = [1,4;2,1;3,2;4,19;5,3;6,7;7,8;7,13;8,9;9,5;9,12;10,22;10,24;...
	%    11,10;11,25;12,11;12,15;13,10;13,12;14,11;14,15;15,2;15,17;...
	%    16,14;17,16;18,16;18,21;19,18;20,14;20,25;21,20;22,6;23,22;...
	%    24,23;24,26;25,24;25,27;26,27;27,28;28,20];
	%%% Roadblocks
	%block = [7,9,14,15,22,24,25];
	%%% Plot
	%figure;
	%hold on;
	%arr = .2;
	%for a = 1:size(A,1) 
	%    D(A(a,1),A(a,2)) = 1;
	%    x0 = lon(indV(A(a,1)));
	%    x1 = lon(indV(A(a,2)));
	%    y0 = lat(indV(A(a,1)));
	%    y1 = lat(indV(A(a,2)));
	%    plot([x0,x1],[y0,y1],'k.-');
	%    plot([arr*x0+(1-arr)*x1,x1],[arr*y0+(1-arr)*y1,y1],'k','LineWidth',2);
	%end
	%for n = 1:28
	%    text(lon(indV(n)),lat(indV(n)),num2str(n));
	%end
	%for b = 1:numel(block)
	%    plot(lon(indV(block(b))),lat(indV(block(b))),'ro');
	%end
	%minx = -7.40135e7;
	%maxx = -7.40085e7;
	%miny = 4.07040e7;
	%maxy = 4.07090e7;
	%axis([minx,maxx,miny,maxy]);
	%set(gca,'XTick',[],'YTick',[]);
	%daspect([1,1,1]);
	%axis off

\begin{figure}
	\begin{tikzpicture}[every node/.style={inner sep=0,outer sep=0},scale=10,->,>=stealth',shorten >=1pt]
		\node [draw,circle,minimum size=4mm] (v1) at (0.0405,0.1340) {\scriptsize $1$};
		\node [draw,circle,minimum size=4mm] (v2) at (0.0362,0.2240) {\scriptsize $2$};
		\node [draw,circle,minimum size=4mm] (v3) at (0,0.3765) {\scriptsize $3$};
		\node [draw,circle,minimum size=4mm] (v4) at (0.0628,-.03) {\scriptsize $4$};
		% \node [draw,circle,minimum size=4mm] (v4) at (0.0628,0.0012) {\scriptsize $4$};
		\node [draw,circle,minimum size=4mm] (v5) at (0.0462,0.4308) {\scriptsize $5$};
		\node [draw,circle,minimum size=4mm] (v6) at (0.4910,1.0000) {\scriptsize $6$};
		\node [draw,circle,minimum size=4mm] (v7) at (0.3983,0.8792) {\scriptsize $7$};
		\node [draw,circle,minimum size=4mm] (v8) at (0.3289,0.7974) {\scriptsize $8$};
		\node [draw,circle,minimum size=4mm] (v9) at (0.2084,0.6493) {\scriptsize $9$};
		\node [draw,circle,minimum size=4mm] (v10) at (0.6095,0.7144) {\scriptsize $10$};
		\node [draw,circle,minimum size=4mm] (v11) at (0.5029,0.5350) {\scriptsize $11$};
		\node [draw,circle,minimum size=4mm] (v12) at (0.3511,0.5941) {\scriptsize $12$};
		\node [draw,circle,minimum size=4mm] (v13) at (0.4809,0.8047) {\scriptsize $13$};
		\node [draw,circle,minimum size=4mm] (v14) at (0.4324,0.2109) {\scriptsize $14$};
		\node [draw,circle,minimum size=4mm] (v15) at (0.1791,0.2157) {\scriptsize $15$};
		\node [draw,circle,minimum size=4mm] (v16) at (0.4320,0.1407) {\scriptsize $16$};
		\node [draw,circle,minimum size=4mm] (v17) at (0.1971,0.1414) {\scriptsize $17$};
		\node [draw,circle,minimum size=4mm] (v18) at (0.4271,0.0313) {\scriptsize $18$};
		\node [draw,circle,minimum size=4mm] (v19) at (0.4130,-.03) {\scriptsize $19$};
		% \node [draw,circle,minimum size=4mm] (v19) at (0.4130,0) {\scriptsize $19$};
		\node [draw,circle,minimum size=4mm] (v20) at (0.7411,0.2451) {\scriptsize $20$};
		\node [draw,circle,minimum size=4mm] (v21) at (0.6762,0.1400) {\scriptsize $21$};
		\node [draw,circle,minimum size=4mm] (v22) at (0.7162,0.8391) {\scriptsize $22$};
		\node [draw,circle,minimum size=4mm] (v23) at (0.9734,0.6664) {\scriptsize $23$};
		\node [draw,circle,minimum size=4mm] (v24) at (0.8669,0.5457) {\scriptsize $24$};
		\node [draw,circle,minimum size=4mm] (v25) at (0.7926,0.4090) {\scriptsize $25$};
		\node [draw,circle,minimum size=4mm] (v26) at (1.0000,0.4474) {\scriptsize $26$};
		\node [draw,circle,minimum size=4mm] (v27) at (0.9444,0.3421) {\scriptsize $27$};
		\node [draw,circle,minimum size=4mm] (v28) at (0.9487,0.2757) {\scriptsize $28$};
		\path [->] (v1) edge node [right] {\scriptsize Broadway} (v4);
		\path [->] (v2) edge node [left] {} (v1);
		\path [->] (v3) edge node [left] {} (v2);
		\path [->] (v4) edge node [below] {\scriptsize Stone} (v19);
		\path [->] (v5) edge node [left] {} (v3);
		\path [->] (v6) edge node [left] {\scriptsize Broadway} (v7);
		\path [->] (v7) edge node [left] {\scriptsize Broadway} (v8);
		\path [->,dashed] (v7) edge node [right] {\scriptsize Wall} (v13);
		\path [->] (v8) edge node [left] {\scriptsize Broadway} (v9);
		\path [->] (v9) edge node[left] {\scriptsize Broadway} (v5);
		\path [->,dashed] (v9) edge node [above] {\scriptsize Exchange} (v12);
		\path [->,dashed] (v10) edge node [right] {\scriptsize Broad} (v22);
		\path [->,dashed] (v10) edge node [right] {\scriptsize Wall} (v24);
		\path [->] (v11) edge node [left] {\scriptsize Broad} (v10);
		\path [->,dashed] (v11) edge node [above] {\scriptsize Exchange} (v25);
		\path [->] (v12) edge node [above] {\scriptsize Exchange} (v11);
		\path [->,dashed] (v12) edge node [left] {\scriptsize New} (v15);
		\path [->] (v13) edge node [right] {\scriptsize Wall} (v10);
		\path [->] (v13) edge node [left] {\scriptsize New} (v12);
		\path [->,dashed] (v14) edge node [left] {\scriptsize Broad} (v11);
		\path [->] (v14) edge node [below] {\scriptsize Beaver} (v15);
		\path [->] (v15) edge node [below] {\scriptsize Beaver} (v2);
		\path [->] (v15) edge node [below] {} (v17);
		\path [->] (v16) edge node [below] {} (v14);
		\path [->] (v17) edge node [below] {\scriptsize Marketfield} (v16);
		\path [->] (v18) edge node [left] {\scriptsize Broad} (v16);
		\path [->] (v18) edge node [above] {\scriptsize William} (v21);
		\path [->] (v19) edge node [below] {} (v18);
		\path [->] (v20) edge node [below] {\scriptsize Beaver} (v14);
		\path [->] (v20) edge node [left] {\scriptsize William} (v25);
		\path [->] (v21) edge node [right] {\scriptsize William} (v20);
		\path [->] (v22) edge node [right] {\scriptsize Pine} (v6);
		\path [->] (v23) edge node [right] {\scriptsize Pine} (v22);
		\path [->] (v24) edge node [left] {\scriptsize William} (v23);
		\path [->] (v24) edge node [right] {\scriptsize Wall} (v26);
		\path [->] (v25) edge node [below] {\scriptsize William} (v24);
		\path [->] (v25) edge node [below] {\scriptsize Exchange} (v27);
		\path [->] (v26) edge node [left] {\scriptsize Hanover} (v27);
		\path [->] (v27) edge node [below] {} (v28);
		\path [->] (v28) edge node [below] {\scriptsize Beaver} (v20);
	\end{tikzpicture}
	\caption{\label{fig:FiDi} Embedding of the directed graph $D$ defined by street directions for the heart of New York City's financial district (vertex indices and embedding coordinates based on \texttt{http://www.dis.uniroma1.it/challenge9/download.shtml}; street directions obtained from Google Maps [except for New Street, for which no direction was easily determined]). After 9/11, roadblocks were emplaced at positions corresponding to the periphery of dashed arcs to create a multi-block security zone surrounding the New York Stock Exchange (see \texttt{http://caselaw.findlaw.com/ny-supreme-court/1037146.html}).}
\end{figure}

\begin{figure}
	\begin{tikzpicture}[every node/.style={inner sep=0,outer sep=0},scale=10,->,>=stealth',shorten >=1pt]
		\node [draw,color=white,circle,minimum size=4mm] (v1) at (0.0405,0.1340) {\scriptsize $1$};
		\node [draw,circle,minimum size=4mm] (v2) at (0.0362,0.2240) {\scriptsize $2$};
		\node [draw,color=white,circle,minimum size=4mm] (v3) at (0,0.3765) {\scriptsize $3$};
		\node [draw,color=white,circle,minimum size=4mm] (v4) at (0.0628,-.03) {\scriptsize $4$};
		% \node [draw,color=white,circle,minimum size=4mm] (v4) at (0.0628,0.0012) {\scriptsize $4$};
		\node [draw,color=white,circle,minimum size=4mm] (v5) at (0.0462,0.4308) {\scriptsize $5$};
		\node [draw,color=white,circle,minimum size=4mm] (v6) at (0.4910,1.0000) {\scriptsize $6$};
		\node [draw,circle,minimum size=4mm] (v7) at (0.3983,0.8792) {\scriptsize $7$};
		\node [draw,color=white,circle,minimum size=4mm] (v8) at (0.3289,0.7974) {\scriptsize $8$};
		\node [draw,circle,minimum size=4mm] (v9) at (0.2084,0.6493) {\scriptsize $9$};
		\node [draw,circle,minimum size=4mm] (v14) at (0.4324,0.2109) {\scriptsize $14$};
		\node [draw,circle,minimum size=4mm] (v15) at (0.1791,0.2157) {\scriptsize $15$};
		\node [draw,circle,minimum size=4mm] (v16) at (0.4320,0.1407) {\scriptsize $16$};
		\node [draw,color=white,circle,minimum size=4mm] (v17) at (0.1971,0.1414) {\scriptsize $17$};
		\node [draw,circle,minimum size=4mm] (v18) at (0.4271,0.0313) {\scriptsize $18$};
		\node [draw,color=white,circle,minimum size=4mm] (v19) at (0.4130,-.03) {\scriptsize $19$};
		% \node [draw,circle,minimum size=4mm] (v19) at (0.4130,0) {\scriptsize $19$};
		\node [draw,color=white,circle,minimum size=4mm] (v20) at (0.7411,0.2451) {\scriptsize $20$};
		\node [draw,color=white,circle,minimum size=4mm] (v21) at (0.6762,0.1400) {\scriptsize $21$};
		\node [draw,circle,minimum size=4mm] (v22) at (0.7162,0.8391) {\scriptsize $22$};
		\node [draw,color=white,circle,minimum size=4mm] (v23) at (0.9734,0.6664) {\scriptsize $23$};
		\node [draw,circle,minimum size=4mm] (v24) at (0.8669,0.5457) {\scriptsize $24$};
		\node [draw,circle,minimum size=4mm] (v25) at (0.7926,0.4090) {\scriptsize $25$};
		\node [draw,color=white,circle,minimum size=4mm] (v26) at (1.0000,0.4474) {\scriptsize $26$};
		\node [draw,color=white,circle,minimum size=4mm] (v27) at (0.9444,0.3421) {\scriptsize $27$};
		\node [draw,circle,minimum size=4mm] (v28) at (0.9487,0.2757) {\scriptsize $28$};
		\foreach \from/\to in {
			v9/v2, v15/v2, v7/v9, v9/v22, v7/v22, v14/v22, v24/v22, v22/v7, v16/v14, v28/v14, v9/v15, v7/v15, v14/v15, v15/v16, v18/v16, v2/v18, v18/v28, v25/v28, v9/v24, v7/v24, v14/v24, v25/v24, v9/v25, v7/v25, v14/v25, v28/v25, v24/v25}
		\draw (\from) to (\to);
	\end{tikzpicture}
	\caption{\label{fig:FiDiAbstraction} $D \upuparrows (\ell,L)$ with $D$ as in figure \ref{fig:FiDi}, $\ell^{\times 28}([28]) = (1,1,2,1,2,3,4,4,2,5,5,5,5,6,7,8,7,9,1,10,9,3,11,11,12,11,12,10)$, and $L = [12] \setminus \{5\}$.}
\end{figure}
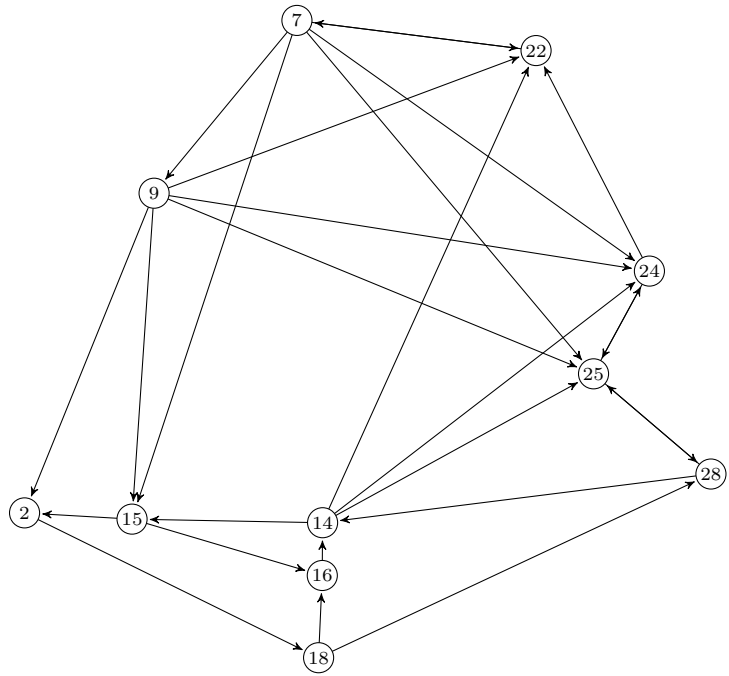

\textsc{Example.} Let $n = 28$ and $p = 0.05$, and suppose that the graph $D$ in figure \ref{fig:FiDi} is a realization of $D_{n,p}$. (Note that $D$ has $|A(D)| = \sum_{x,y} \mu_D(x,y) = 40$ arcs and hence $\frac{|A(D)|}{n(n-1)} \approx 0.0529 \approx p$. Furthermore, 12 vertices have indegree 2 and the rest have indegree 1; similarly, 12 vertices have outdegree 2 and the rest have outdegree 1, so $D_{n,p}$ is at least superficially appropriate here as a model for $D$.) The partial partition $\pi$ corresponding to the path abstraction in figure \ref{fig:FiDiAbstraction} has $n - |\text{supp } \pi| = 4$, $f^4(p) \approx 0.0578$, and $(|\pi^{(1)}|,\dots,|\pi^{(|\pi|)}|) = (4,3,2,2,1,2,1,2,2,3,2)$. The expected number of arcs in $D_{n,p} \upuparrows \pi$ turns out to be approximately 25.9635, whereas the actual number of arcs in $D \upuparrows \pi$ is 27. (Taking instead the value $p = 0.0529$ yields an expected number of arcs approximately equal to 27.4466.) $\Box$

\

\begin{figure}[htbp]
\includegraphics[trim = 0mm 0mm 0mm 0mm, clip, width=90mm,keepaspectratio]{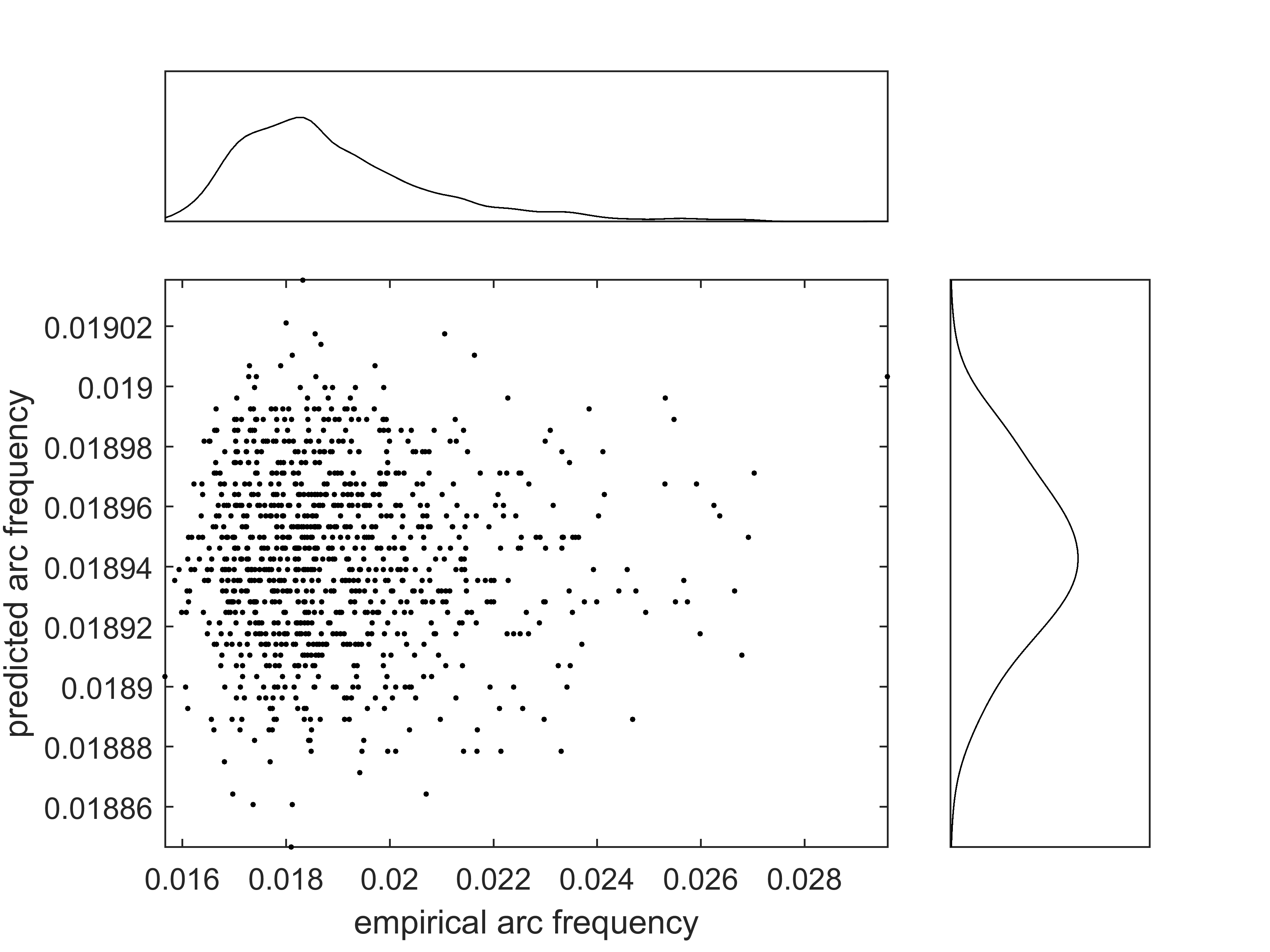}% Here is how to import pix
\caption{ \label{fig:ProbabilityProjections} (Lower left) Scatterplot of empirical arc frequencies and the predictions of \eqref{eq:RandomContraction2} in 1000 realizations of $D_{n,p} \upuparrows U$ with $n = 1000$, $p = 0.01$, and $U$ a uniformly random subset of $[n]$ subject to $|U| = 50$. (Other panels) Marginal kernel density estimates. The maximum of the distribution in the top left panel (obtained with a bandwidth of $2.9649 \times 10^{-4}$) is at $0.01834$; the maximum of the distribution in the lower right panel (obtained with a bandwidth of $1.1114 \times 10^{-5}$) is at $0.01894$. Note that these two maxima are separated by just a few multiples of the larger bandwidth.} 
\end{figure}

\textsc{Example.} Let $n = 1000$, $p = 0.01$, and consider $U$ a uniformly random subset of $[n]$ subject to $|U| = 50$. Figure \ref{fig:ProbabilityProjections} demonstrates \eqref{eq:RandomContraction2} using the estimate $\hat p = \frac{|A(D_{n,p})|}{n(n-1)}$ for 1000 realizations of $D_{n,p} \upuparrows U$. $\Box$

\

    % rng('default'); % for reproducibility
    % n = 1000;
    % p = .01;
    % F = zeros(1000,4);
    % for i = 1:size(F,1)
    %    %%
    %    Dnp = spdiags(sparse(n,1),0,spones(sprand(n,n,p)));
    %    % For more precise results (doesn't make much difference), can replace above with
    %    % Dnp = (rand(n)<p); 
    %    % Dnp = Dnp-diag(diag(Dnp));
    %    N = 50;
    %    U = randperm(n,N);
    %    DnpU = graphdetour(Dnp,U);
    %    %%
    %    num_arcs = sum(Dnp(:));
    %    p_emp = num_arcs/(n*(n-1));
    %    %%
    %    num_arcs_U = sum(DnpU(:));
    %    p_emp_U = num_arcs_U/((n-N)*(n-N-1));
    %    %%
    %    f = p_emp;
    %    for j = 2:N
    %        f = f^2+(1-f^2)*f;
    %    end
    %    F(i,:) = full([p,p_emp,p_emp_U,f]);
    % end
    % tde3 = tde1d(F(:,3),1);
    % tde4 = tde1d(F(:,4),1);
    % %%
    % figure;
    % k = 4;
    % subplot(k,k,1:(k-1))
    % plot(tde3.x,tde3.y,'k');
    % set(gca,'XTick',[],'YTick',[]);
    % xlim([min(tde3.x),max(tde3.x)]);
    % subplot(k,k,setdiff(k:k^2,k*(1:k)))
    % plot(F(:,3),F(:,4),'k.');
    % xlim([min(tde3.x),max(tde3.x)]);
    % ylim([min(tde4.x),max(tde4.x)]);
    % xlabel('empirical arc frequency');
    % ylabel('predicted arc frequency');
    % subplot(k,k,k*(2:k))
    % plot(tde4.y,tde4.x,'k');
    % set(gca,'XTick',[],'YTick',[]);
    % ylim([min(tde4.x),max(tde4.x)]);

\begin{figure}[htbp]
\includegraphics[trim = 0mm 0mm 0mm 0mm, clip, width=90mm,keepaspectratio]{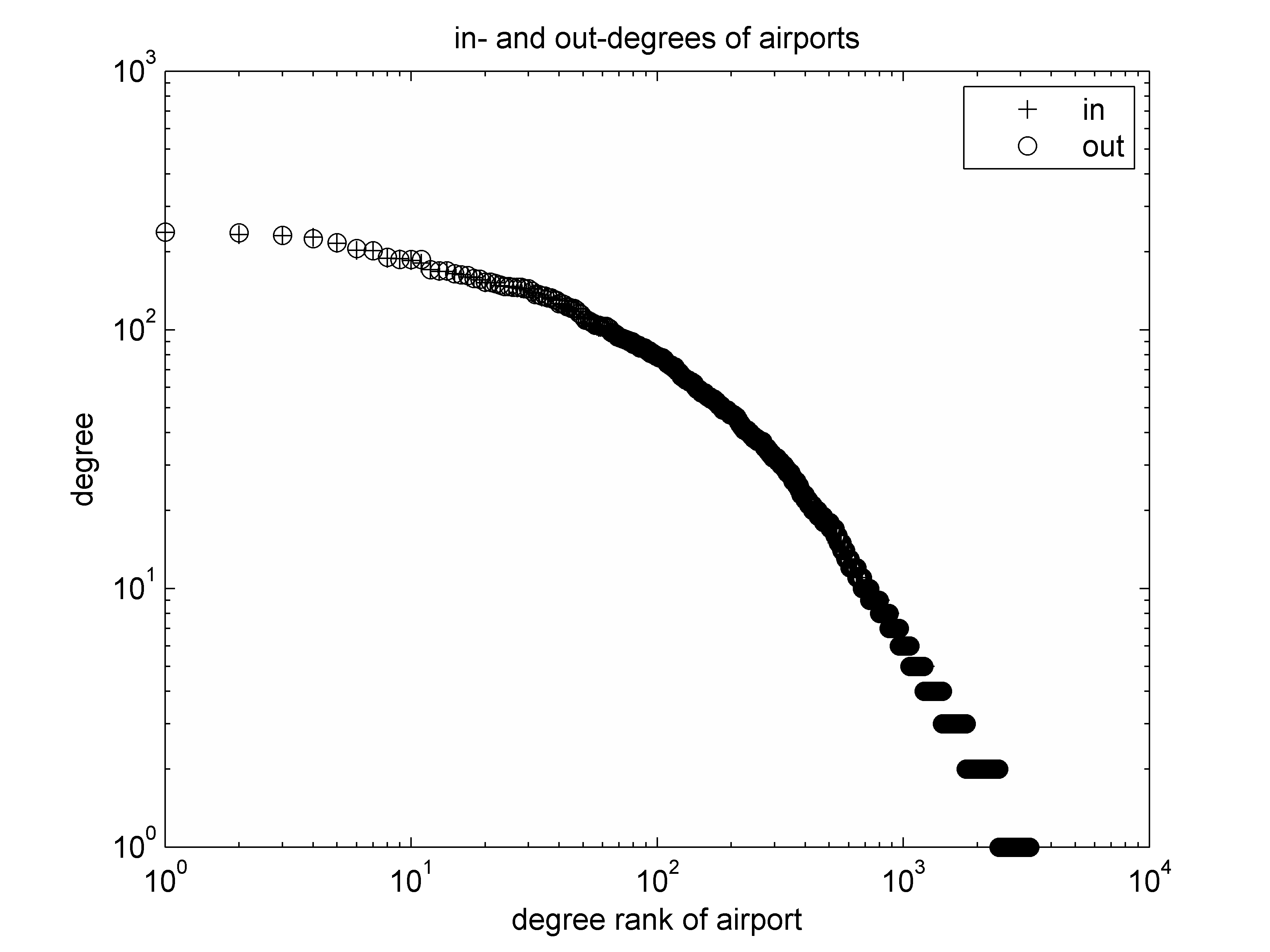}% Here is how to import pix
\includegraphics[trim = 0mm 0mm 0mm 0mm, clip, width=90mm,keepaspectratio]{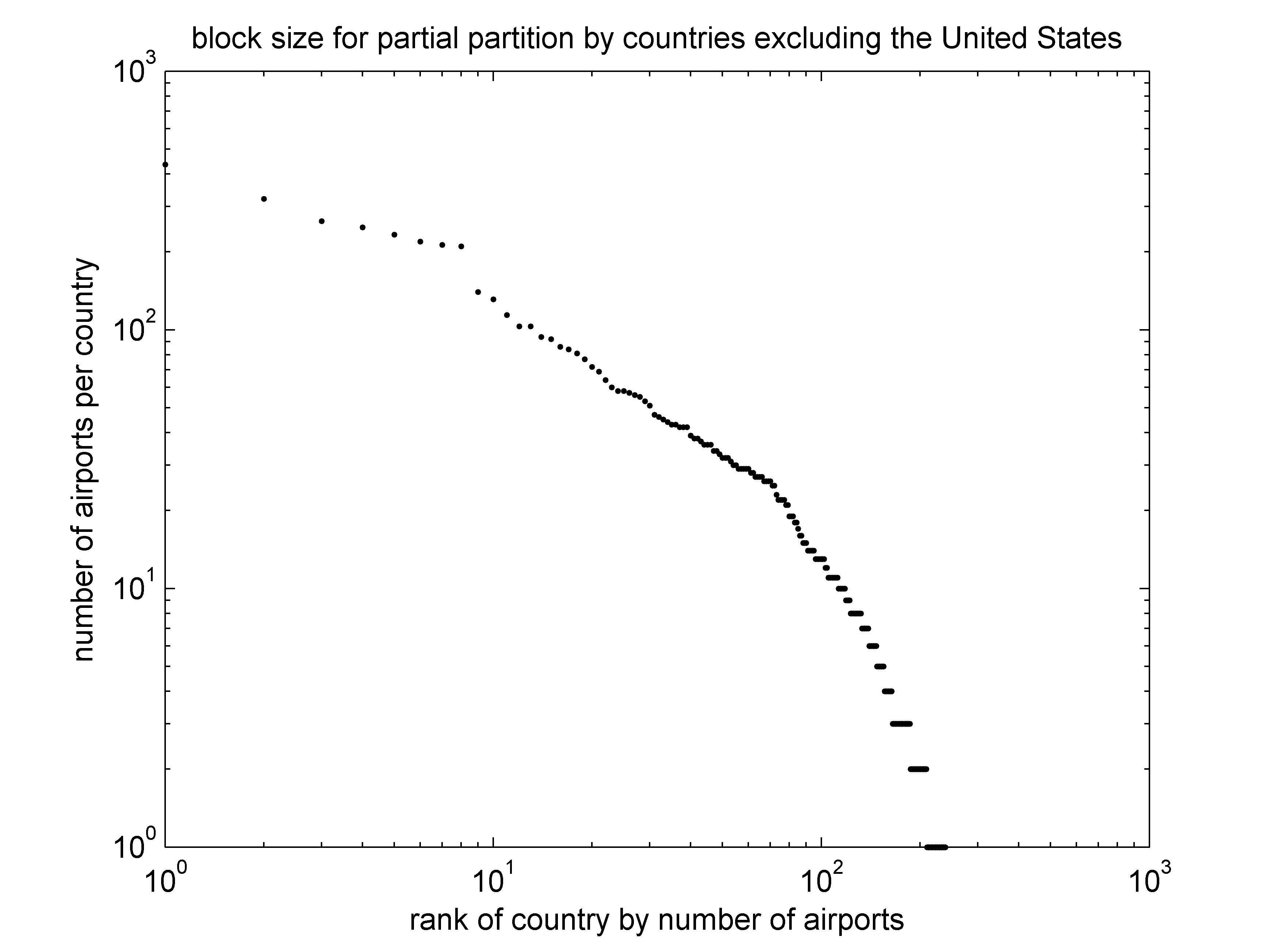}% Here is how to import pix
\caption{ \label{fig:Airports}(L)  In- and out-degrees of airports are very far from uniform. (R) The partial partition $\pi'$ defined by the country in which an airport resides and omitting the United States has very nonuniform block sizes.} 
\end{figure}

\textsc{Example.} As an example for which $D_{n,p}$ is a manifestly inappropriate model, consider the digraph $D'$ with vertices corresponding to airports and arcs corresponding to regularly scheduled passenger flights. We constructed this digraph using data accessed from \texttt{http://openflights.org/data.html} on 3 May 2016, yielding $n' := |V(D')| = 8107$ and $|A(D')| = 37187$ arcs, corresponding to an empirical arc probability of $p' := \frac{|A(D)|}{n(n-1)} \approx 5.6588 \cdot 10^{-4}$. We define $\ell' : V(D') \rightarrow \Lambda'$ to be the map coloring airports by country, so that $|\Lambda'| = 240$. The in- and out-degrees are very far from uniform, as illustrated in figure \ref{fig:Airports}.
Evidently $D_{n',p'}$ is a very poor model for $D'$. If $\pi'$ is the partial partition corresponding to country colors excluding the United States, then $|A(D' \upuparrows \pi')| = 11036$, whereas $\mathbb{E}(|A(D_{n',p'} \upuparrows \pi')|) \approx 5125$. $\Box$

\

In order to get better approximations in such situations it would be necessary to consider a more general random digraph, e.g. the one introduced in \cite{BloznelisGotzeJaworski} (ignoring loops). However, the proof of the preceding theorem exploited the commutativity of contracting and bypassing vertices in an essential way that does not generalize to the random digraph of \cite{BloznelisGotzeJaworski}. For this reason, obtaining a suitable generalization of the theorem appears to require substantially more effort. Furthermore, applying the resulting theorem would appear to require the same sort of set-theoretic operations as actually constructing the path abstraction outright, largely negating its utility as a predictive tool.

\section{\label{sec:TemporalNetworks}Temporal networks}

Digraphs admit a natural temporal generalization called \emph{directed temporal contact networks (DTCNs)}. \cite{MasudaLambiotte} A DTCN with vertex set $V \equiv [n]$ is a pair $(\mathcal{D},\delta)$ where $\mathcal{D}$ is a finite nonempty set and $\delta : \mathcal{D} \rightarrow (V^2 \setminus \Delta(V)) \times \mathbb{R}$ is injective. The \emph{source}, \emph{target}, and \emph{time} maps (respectively denoted $s$, $t$, and $\tau$) are defined so that the following diagram commutes:
\begin{center}
	\begin{tikzpicture}[->,>=stealth',shorten >=1pt]
		\node (v01) at (0,0) {$\mathcal{D}$};
		\node (v02) at (0,-4) {$V^2$};
		\node (v03) at (2,-2) {$(V^2 \setminus \Delta(V)) \times \mathbb{R}$};
		\node (v04) at (4,0) {$\mathbb{R}$};
		\node (v05) at (4,-4) {$V^2 \times \mathbb{R}$};
		\path [->] (v01) edge node [left] {$s \times t$} (v02);
		\path [->] (v01) edge node [above] {$\delta$} (v03);
		\path [->] (v01) edge node [above] {$\tau$} (v04);
		\path [right hook->] (v03) edge node [above] {$i$} (v05);
		\path [->] (v05) edge node [below] {$(\pi_1 \circ \pi_1) \times (\pi_2 \circ \pi_1)$} (v02);
		\path [->] (v05) edge node [right] {$\pi_2$} (v04);
	\end{tikzpicture}
\end{center}
That is, each \emph{contact} $c \in \mathcal{D}$ corresponds to a unique triple $(s(c),t(c),\tau(c))$, and when convenient we identify contacts and their corresponding triples. We may economically indicate the DTCN $(\mathcal{D},\delta)$ merely as $\mathcal{D}$ or $\delta$; context should suffice to remove any potential for ambiguity here. There is an obvious notion of a temporally coherent path which we do not bother to write out formally. 

A naive attempt to generalize the constructions of \S \ref{sec:DetoursBypassesAbstractions} to DTCNs might define $\mathcal{D} \uparrow v \equiv \mathcal{D} \upuparrows v$ as
\begin{equation}
\label{eq:WrongTemporalDetour}
\{c \in \mathcal{D}: s(c) \ne v \ne t(c)\} \cup \{ (j,k,\tau_{vk}) : (j,v,\tau_{jv}), (v,k,\tau_{vk}) \in \mathcal{D} \land \tau_{jv} \le \tau_{vk} \land j \ne k \}.
\end{equation}
However, such a definition yields undesirable behavior, as we show in the following 

\

\begin{figure}
	\begin{tikzpicture}[->,>=stealth',shorten >=1pt]
		\node [draw,circle,minimum size=7mm] (v10) at (-1,-1) {1};
		\node [draw,circle,minimum size=7mm] (v20) at (-1,-2) {2};
		\node [draw,circle,minimum size=7mm] (v30) at (-1,-3) {3};
		\node [draw,circle,minimum size=7mm] (v40) at (-1,-4) {4};
		\node [draw,circle,minimum size=7mm] (v50) at (-1,-5) {5};
		\node at (0,-5.5) {$-\infty$};
		\node at (1,-5.5) {$\tau_1$};
		\node at (2,-5.5) {$\tau_2$};
		\node at (3,-5.5) {$\tau_3$};
		\node at (4,-5.5) {$\tau_4$};
		\node at (5,-5.5) {$\infty$};
		\coordinate (v1a) at (0,-1);
		\coordinate (v2a) at (0,-2);
		\coordinate (v3a) at (0,-3);
		\coordinate (v4a) at (0,-4);
		\coordinate (v5a) at (0,-5);
		\coordinate (v11) at (1,-1);
		\coordinate (v21) at (1,-2);
		\coordinate (v31) at (1,-3);
		\coordinate (v41) at (1,-4);
		\coordinate (v51) at (1,-5);
		\coordinate (v42) at (2,-4);
		\coordinate (v52) at (2,-5);
		\coordinate (v23) at (3,-2);
		\coordinate (v53) at (3,-5);
		\coordinate (v14) at (4,-1);
		\coordinate (v24) at (4,-2);
		\coordinate (v34) at (4,-3);
		\coordinate (v44) at (4,-4);
		\coordinate (v54) at (4,-5);
		\coordinate (v1z) at (5,-1);
		\coordinate (v2z) at (5,-2);
		\coordinate (v3z) at (5,-3);
		\coordinate (v4z) at (5,-4);
		\coordinate (v5z) at (5,-5);
		\foreach \from/\to in {
			v1a/v11, v11/v1z, v2a/v23, v23/v2z, v3a/v34, v34/v3z, v4a/v41, v41/v42, v42/v44, v44/v4z, v5a/v52, v52/v53, v53/v5z}
		\draw (\from) to (\to);
		\draw (v11) [out=-60,in=60,looseness=1] to (v41);
		\draw (v52) [out=60,in=-60,looseness=1] to (v42);
		\draw (v23) [out=-60,in=60,looseness=1] to (v53);
		\draw (v44) [out=60,in=-60,looseness=1] to (v34);
	\end{tikzpicture}
	\caption{\label{fig:DTCNtemporalgraph} Temporal digraph of the DTCN $\mathcal{D} := \{(1,4,\tau_1), (5,4,\tau_2), (2,5,\tau_3), (4,3,\tau_4)\}$ with $\tau_1 < \tau_2 < \tau_3 < \tau_4$. Note that there is a temporally coherent path from 1 to 3, but not from 2 to 3.}
\end{figure}
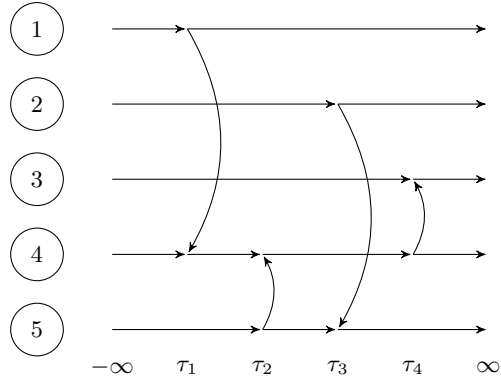

\textsc{Example.} Consider the DTCN $\mathcal{D} := \{(1,4,\tau_1), (5,4,\tau_2), (2,5,\tau_3), (4,3,\tau_4)\}$ with $\tau_1 < \tau_2 < \tau_3 < \tau_4$. Using the naive definition \eqref{eq:WrongTemporalDetour} for $\mathcal{D} \uparrow v$ leads to $(\mathcal{D} \uparrow 4) \uparrow 5 = \{(1,3,\tau_4),(2,3,\tau_4)\} \ne (\mathcal{D} \uparrow 5) \uparrow 4 = \{(1,3,\tau_4)\}$. As figure \ref{fig:DTCNtemporalgraph} illustrates, only the latter corresponds to the desirable result for $\mathcal{D} \uparrow \{4,5\}$. $\Box$

\

An approach that manifestly yields the desired construction deals with the \emph{temporal digraph} of $\mathcal{D}$ (see figure \ref{fig:DTCNtemporalgraph} for an example), defined as the digraph $T(\mathcal{D})$ with vertex and arc sets
\begin{eqnarray}
\label{eq:TemporalDigraphVertices}
V(T(\mathcal{D})) & := & \{(v, -\infty) : v \in V\} \cup \{(v,\tau(c)) : (v,c) \in V \times \mathcal{D} \land (s(c) = v \lor t(c) = v) \} \cup \{(v, \infty) : v \in V\} \\
\label{eq:TemporalDigraphArcs}
A(T(\mathcal{D})) & := & \{((s(c),\tau(c)),(t(c),\tau(c))) : c \in \mathcal{D} \} \cup \{ ((v,\tau_{j-1}^{@v}),(v,\tau_j^{@v})) : v \in V, j \in [|\mathcal{D}@v|-1] \}
\end{eqnarray} 
where the \emph{temporal fiber at $v$ is}
\begin{equation}
\label{eq:TemporalFiber}
\mathcal{D}@v := \{-\infty\} \cup \{\tau(c) : c \in \mathcal{D} \land s(c) = v \} \cup \{\tau(c) : c \in \mathcal{D} \land t(c) = v \} \cup \{\infty\} \equiv \{\tau_j^{@v}\}_{j = 0}^{|\mathcal{D}@v|-1}. \nonumber
\end{equation}
Note that $|V(T(\mathcal{D}))| = \sum_v |\mathcal{D}@v| \le 2|V| + 2|\mathcal{D}|$ and $|A(T(\mathcal{D}))| = |V(T(\mathcal{D}))| - |V| + |\mathcal{D}| \le |V| + 3 |\mathcal{D}|$, so that the temporal digraph of a DTCN can be formed with only linear overhead. 

Call the two sets in the union on the RHS of \eqref{eq:TemporalDigraphArcs} the \emph{spatial} and \emph{temporal arcs} of $T(\mathcal{D})$, respectively. Now for $U \subseteq V$, consider
\begin{equation}
T(\mathcal{D}) \upuparrows \bigcup_{u \in U} \left ( \{u\} \times \mathcal{D}@u \right ). \nonumber
\end{equation}
Each of the non-temporal arcs in this digraph is of the form $((v,\tau_j^{@v}),(v',\tau_{j'}^{@v'}))$ for $v \ne v'$ and corresponds to a triple $(v,v',\tau_j^{@v} \lor \tau_{j'}^{@v'})$. The set of all such triples defines $\mathcal{D} \uparrow U \equiv \mathcal{D} \upuparrows U$. That is,
\begin{equation}
\label{eq:TemporalDetour}
\mathcal{D} \uparrow U \equiv \mathcal{D} \upuparrows U := \left \{ (v,v',\tau_j^{@v} \lor \tau_{j'}^{@v'}) : v \ne v' \land ((v,\tau_j^{@v}),(v',\tau_{j'}^{@v'})) \in A \left ( T(\mathcal{D}) \upuparrows \bigcup_{u \in U} \left ( \{u\} \times \mathcal{D}@u \right ) \right ) \right \}.
\end{equation}

\

\textsc{Proposition}. If $\tau$ is a constant map, then $\mathcal{D} \uparrow U$ can be identified with $D \uparrow U$, where here $D$ indicates the obvious digraph corresponding to $\mathcal{D}$. $\Box$

\

	%n = 5; 
	%dtcn = [1,4,1;5,4,2;2,5,3;4,3,4]
	%U = [4,5];
	%dtcndetour(n,dtcn,U(1)) % = [1,3,4;2,5,3;5,3,4]
	%dtcndetour(n,dtcn,U(2)) % = [1,4,1;4,3,4]
	%dtcndetour(n,dtcndetour(n,dtcn,U(1)),U(2)) % = [1,3,4;2,3,4]
	%dtcndetour(n,dtcndetour(n,dtcn,U(2)),U(1)) % = [1,3,4]
	%dtcndetour(n,dtcn,U) % = [1,3,4]

\textsc{Example.} Again, consider the DTCN $\mathcal{D} := \{(1,4,\tau_1), (5,4,\tau_2), (2,5,\tau_3), (4,3,\tau_4)\}$ with $\tau_1 < \tau_2 < \tau_3 < \tau_4$. Using \eqref{eq:TemporalDetour} yields $\mathcal{D} \uparrow \{4,5\} = \{(1,3,\tau_4)\}$, as desired. However, it is still the case that $(\mathcal{D} \uparrow \{4\}) \uparrow \{5\} = \{(1,3,\tau_4),(2,3,\tau_4)\} \ne (\mathcal{D} \uparrow \{5\}) \uparrow \{4\} = \{(1,3,\tau_4)\}$, just as before. $\Box$

\

The preceding examples show that although \eqref{eq:TemporalDetour} is certainly a reasonable definition for $\mathcal{D} \uparrow U$, any reasonable definition of detours/bypasses for DTCNs will lead to noncommutativity that is not present for digraphs. However, there is still a well-defined notion of path abstraction for DTCNs (which necessarily will not commute with successive detours/bypasses) due to the following

\

\textsc{Lemma.} Detours/bypasses commute with vertex contractions for DTCNs.

\

\textsc{Sketch of proof.} Let $U \cap \{v,w\} = \varnothing$ and $U \cup \{v,w\} \subseteq V$. The vertex contraction $\mathcal{D} / \{v,w\}$ is defined in the obvious way: triples of the form $(v,x,\tau_{vx})$ and $(x,v,\tau_{xv})$ for $x \not \in \{v,w\}$ are replaced with $(\{v,w\},x,\tau_{vx})$ and $(x,\{v,w\},\tau_{xv})$, respectively, and similarly for triples involving $w$. Thus $(\mathcal{D} / \{v,w\})@\{v,w\} = (\mathcal{D}@v) \cup (\mathcal{D}@w)$, so additional vertices and temporal arcs are generated in the formation of $T(\mathcal{D} / \{v,w\})$. However, in the formation of $T(\mathcal{D})$, replacing both of the  temporal fibers $\mathcal{D}@v$ and $\mathcal{D}@w$ with $(\mathcal{D} / \{v,w\})@\{v,w\}$ has no material effect on the subsequent formation of $\mathcal{D} \uparrow U$. The lemma now reduces to the already established version for digraphs by contracting vertices with the same time coordinate in each copy of $(\mathcal{D} / \{v,w\})@\{v,w\}$. $\Box$

\

The surprising noncommutativity of detours/bypasses for DTCNs is not the only difference from the situation for digraphs.

\

\textsc{Example.} There are at least two random DTCNs that are obvious analogues of $D_{n,p}$ (cf. \S \ref{sec:RandomDigraphs}): 
\begin{itemize}
\item $\mathcal{D}_{n,p}^{(u)}$, with sources and targets corresponding to $D_{n,p}$ and times uniformly random in $[0,1]$;
\item $\mathcal{D}_{n,p}^{(P)}$, with contacts between $x \ne y$ Poisson distributed over $[0,1]$ with rate $p$. That is, the probability of a contact between $x \ne y$ in an interval of infinitesmal duration $d \tau$ is given by $p \cdot d \tau$.
\end{itemize}
It is easy to see that both of these have an expected number of contacts equal to $p \cdot n (n-1)$. Furthermore, for the regime of interest $p \ll 1$, these two random DTCNs can be expected to behave quite similarly, akin to the Erd\H{o}s-R\'enyi and Gilbert random graphs.

\begin{figure}[htbp]
\includegraphics[trim = 0mm 0mm 0mm 0mm, clip, width=90mm,keepaspectratio]{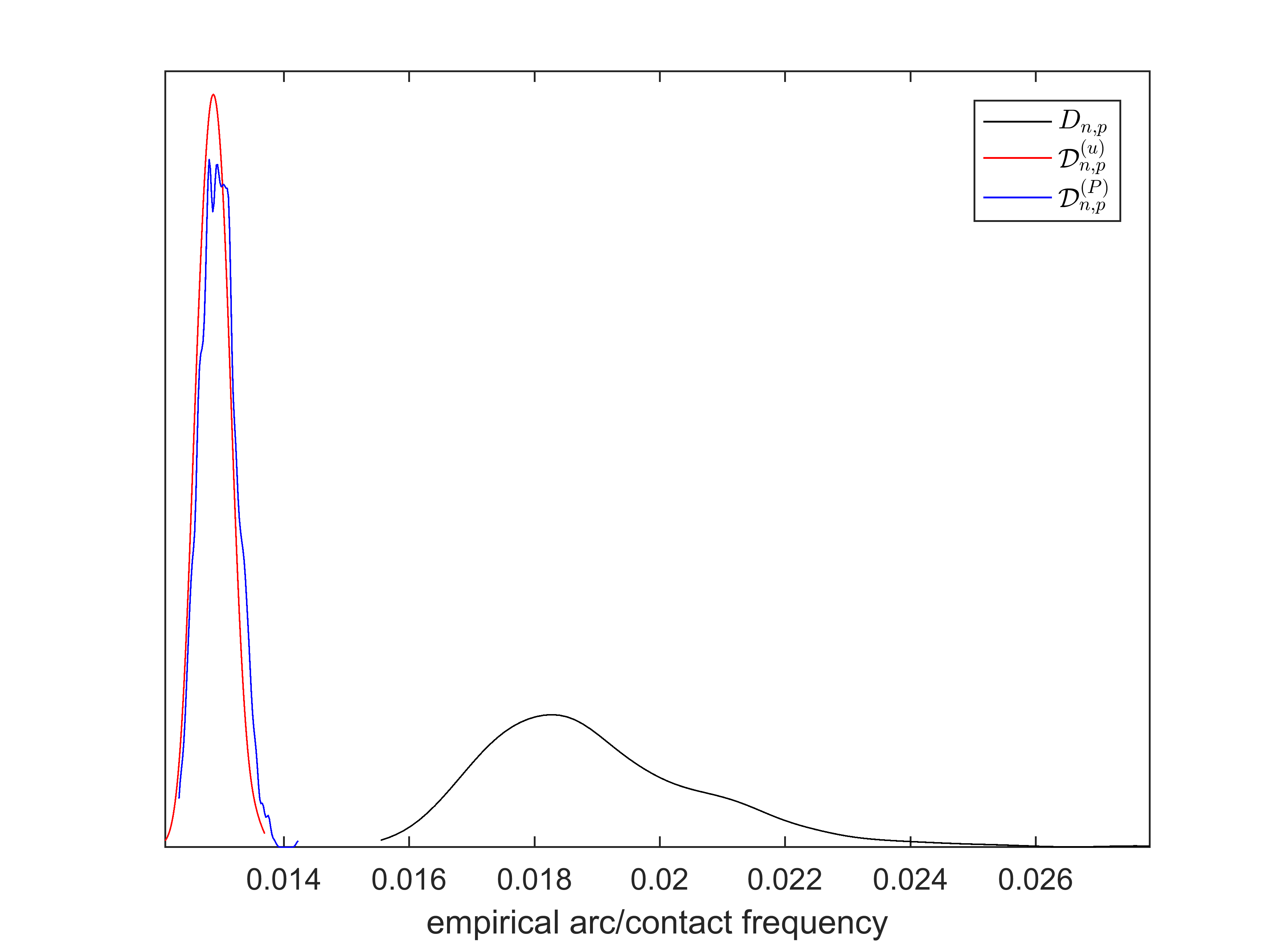}% Here is how to import pix
\caption{ \label{fig:DTCNProbability} Kernel density estimates of empirical arc/contact frequencies for 1000 realizations of $D_{n,p} \upuparrows U$ (cf. figure \ref{fig:ProbabilityProjections}) and $\mathcal{D}_{n,p}^{(\cdot)} \uparrow U$ with $n = 1000$, $p = 0.01$, and $U$ a uniformly random subset of $[n]$ subject to $|U| = 50$.} 
\end{figure}

Rather than attempting to develop analytical results, we proceed directly to numerics. The basic observation from figure \ref{fig:DTCNProbability} is that the number of contacts in $\mathcal{D}_{n,p}^{(\cdot)} \uparrow U$ is much less than the number of arcs in $D_{n,p} \uparrow U$, because there are fewer temporally coherent paths between two vertices in $\mathcal{D}_{n,p}^{(\cdot)}$ than there are ordinary paths between the same two vertices in $D_{n,p}$. In particular, given $\{x_j\}_{j = 0}^\ell$, the probability that the path $x_0 \rightarrow \dots \rightarrow x_\ell$ exists in $D_{n,p}$ is $p^\ell$, whereas the probability that a temporally coherent version of the same path exists in $\mathcal{D}_{n,p}^{(\cdot)}$ is $p^\ell/\ell!$
$\Box$

\acknowledgements

The author is grateful to Mukesh Dalal for proposing the idea of the paper, and to Yingbo Song for patiently allowing it to unfold. This material is based upon work supported by the Defense Advanced Research Projects Agency (DARPA) and the Air Force Research Laboratory (AFRL). Any opinions, findings and conclusions or recommendations expressed in this material are those of the author(s) and do not necessarily reflect the views of DARPA or AFRL.

\appendix

\section{\label{sec:CommutativeDetours}Case analysis proof of $(D \uparrow v) \uparrow w = (D \uparrow w) \uparrow v$ for digraphs}

There are five cases: 1) $w = v$; 2) $w \in V(D)_v^-$; 3) $w \in V(D)_v^\pm$; 4) $w \in V(D)_v^+$, and 5) $w \in V(D)_v^0$. Case 1) is trivial, but it is still worth observing that $V(D \uparrow v)_v^0 = V(D) \setminus \{v\}$, whereupon we formally obtain the intuitively obvious fact that $(D \uparrow v) \uparrow v = D \uparrow v$. Note that case 2) is equivalent to $v \in V(D)_w^+$, so that cases 2) and 4) are equivalent by symmetry; we will address the latter.

Before proceeding with the remaining cases 3), 4), and 5), let us first define 
\begin{equation}
Z(D)_v := [V(D) \times \{v\}] \cup [\{v\} \times V(D)] \nonumber
\end{equation}
and
\begin{equation}
U(D)_v := P(D)_v \times S(D)_v. \nonumber
\end{equation}
By construction, we have that $\mu_{D \uparrow v}(Z_v) \equiv 0$ and $\mu_{D \uparrow v}(U_v \setminus \Delta(V(D))) \equiv 1$. Now
\begin{equation}
Z(D \uparrow v)_w \cup Z(D)_v = [V(D) \times \{v,w\}] \cup [\{v,w\} \times V(D)] = Z(D \uparrow w)_v \cup Z(D)_w \nonumber
\end{equation}
so it suffices to show that
\begin{equation}
U(D \uparrow v)_w \cup [U(D)_v \setminus Z(D \uparrow v)_w] = U(D \uparrow w)_v \cup [U(D)_w \setminus Z(D \uparrow w)_v] \nonumber
\end{equation}
or equivalently (writing $\triangle$ as usual for symmetric difference)
\begin{equation}
\label{eq:detourU}
[U(D \uparrow v)_w \cup U(D)_v] \ \triangle \ [U(D \uparrow w)_v \cup U(D)_w] \subseteq [V(D) \times \{v,w\}] \cup [\{v,w\} \times V(D)].
\end{equation}

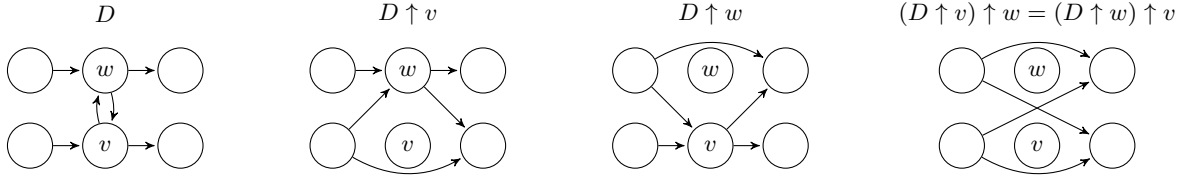
\begin{figure}
	\begin{tikzpicture}[->,>=stealth',shorten >=1pt]
		\node (D) at (1,1.75) {$D$};
		\node [draw,circle,minimum size=6mm] (v1) at (0,0) {};
		\node [draw,circle,minimum size=6mm] (v2) at (0,1) {};
		\node [draw,circle,minimum size=6mm] (v3) at (1,0) {$v$};
		\node [draw,circle,minimum size=6mm] (v4) at (1,1) {$w$};
		\node [draw,circle,minimum size=6mm] (v5) at (2,0) {};
		\node [draw,circle,minimum size=6mm] (v6) at (2,1) {};
		\foreach \from/\to in {
			v1/v3, v2/v4, v3/v5, v4/v6}
		\draw (\from) to (\to);
		\draw (v3) [out=105,in=-105,looseness=1] to (v4);
		\draw (v4) [out=-75,in=75,looseness=1] to (v3);
		\draw (v1) [color=white,out=-30,in=-150,looseness=1] to (v5);	% vertical alignment hack
	\end{tikzpicture}
	\quad \quad \quad \quad
	\begin{tikzpicture}[->,>=stealth',shorten >=1pt]
		\node (D) at (1,1.75) {$D \uparrow v$};
		\node [draw,circle,minimum size=6mm] (v1) at (0,0) {};
		\node [draw,circle,minimum size=6mm] (v2) at (0,1) {};
		\node [draw,circle,minimum size=6mm] (v3) at (1,0) {$v$};
		\node [draw,circle,minimum size=6mm] (v4) at (1,1) {$w$};
		\node [draw,circle,minimum size=6mm] (v5) at (2,0) {};
		\node [draw,circle,minimum size=6mm] (v6) at (2,1) {};
		\foreach \from/\to in {
			v1/v4, v2/v4, v4/v5, v4/v6}
		\draw (\from) to (\to);
		\draw (v1) [out=-30,in=-150,looseness=1] to (v5);
	\end{tikzpicture}
	\quad \quad \quad \quad
	\begin{tikzpicture}[->,>=stealth',shorten >=1pt]
		\node (D) at (1,1.75) {$D \uparrow w$};
		\node [draw,circle,minimum size=6mm] (v1) at (0,0) {};
		\node [draw,circle,minimum size=6mm] (v2) at (0,1) {};
		\node [draw,circle,minimum size=6mm] (v3) at (1,0) {$v$};
		\node [draw,circle,minimum size=6mm] (v4) at (1,1) {$w$};
		\node [draw,circle,minimum size=6mm] (v5) at (2,0) {};
		\node [draw,circle,minimum size=6mm] (v6) at (2,1) {};
		\foreach \from/\to in {
			v1/v3, v2/v3, v3/v5, v3/v6}
		\draw (\from) to (\to);
		\draw (v2) [out=30,in=150,looseness=1] to (v6);
		\draw (v1) [color=white,out=-30,in=-150,looseness=1] to (v5);	% vertical alignment hack
	\end{tikzpicture}
	\quad \quad \quad
	\begin{tikzpicture}[->,>=stealth',shorten >=1pt]
		\node (D) at (1,1.75) {$(D \uparrow v) \uparrow w = (D \uparrow w) \uparrow v$};
		\node [draw,circle,minimum size=6mm] (v1) at (0,0) {};
		\node [draw,circle,minimum size=6mm] (v2) at (0,1) {};
		\node [draw,circle,minimum size=6mm] (v3) at (1,0) {$v$};
		\node [draw,circle,minimum size=6mm] (v4) at (1,1) {$w$};
		\node [draw,circle,minimum size=6mm] (v5) at (2,0) {};
		\node [draw,circle,minimum size=6mm] (v6) at (2,1) {};
		\foreach \from/\to in {
			v1/v6, v2/v5}
		\draw (\from) to (\to);
		\draw (v2) [out=30,in=150,looseness=1] to (v6);
		\draw (v1) [out=-30,in=-150,looseness=1] to (v5);
	\end{tikzpicture}
	\caption{\label{fig:Case3Cartoon} Cartoon for case 3 of the lemma.} 
\end{figure}

\emph{Case 3: $w \in V(D)_v^\pm$.} In this case (see figure \ref{fig:Case3Cartoon} for a cartoon and note that) we have the following identities:
\begin{eqnarray}
P(D \uparrow v)_w & = & P(D)_v \cup [P(D)_w \setminus \{v\}]; \nonumber \\
S(D \uparrow v)_w & = & S(D)_v \cup [S(D)_w \setminus \{v\}]; \nonumber \\
P(D \uparrow w)_v & = & [P(D)_v \setminus \{w\}] \cup P(D)_w; \nonumber \\
S(D \uparrow w)_v & = & [S(D)_v \setminus \{w\}] \cup S(D)_w. \nonumber
\end{eqnarray}
From these it follows that 
\begin{eqnarray}
U(D \uparrow v)_w \cup U(D)_v & \equiv & [P(D \uparrow v)_w \times S(D \uparrow v)_w] \cup [P(D)_v \times S(D)_v] \nonumber \\
& = & [P(D)_v \times S(D)_v] \cup [P(D)_v \times (S(D)_w \setminus \{v\})] \nonumber \\
& \ & \cup \ [(P(D)_w \setminus \{v\}) \times S(D)_v] \cup [(P(D)_w \setminus \{v\}) \times (S(D)_w \setminus \{v\})] \nonumber
\end{eqnarray}
and by symmetry
\begin{eqnarray}
U(D \uparrow w)_v \cup U(D)_w & = & [P(D)_w \times S(D)_w] \cup [P(D)_w \times (S(D)_v \setminus \{w\})] \nonumber \\
& \ & \cup \ [(P(D)_v \setminus \{w\}) \times S(D)_w] \cup [(P(D)_v \setminus \{w\}) \times (S(D)_v \setminus \{w\})] \nonumber
\end{eqnarray}
so upon inspection \eqref{eq:detourU} is satisfied and case 3) is done.

\emph{Case 4: $w \in V(D)_v^+$.} In this case (see figure \ref{fig:Case4Cartoon} for a cartoon and note that) we have the following identities:
\begin{eqnarray}
P(D \uparrow v)_w & = & P(D)_v \cup [P(D)_w \setminus \{v\}]; \nonumber \\
S(D \uparrow v)_w & = & S(D)_w; \nonumber \\
P(D \uparrow w)_v & = & P(D)_v; \nonumber \\
S(D \uparrow w)_v & = & [S(D)_v \setminus \{w\}] \cup S(D)_w. \nonumber
\end{eqnarray}
From these it follows that 
\begin{eqnarray}
U(D \uparrow v)_w \cup U(D)_v & \equiv & [P(D \uparrow v)_w \times S(D \uparrow v)_w] \cup [P(D)_v \times S(D)_v] \nonumber \\
& = & [P(D)_v \times S(D)_w] \cup [(P(D)_w \setminus \{v\}) \times S(D)_w] \cup [P(D)_v \times S(D)_v] \nonumber
\end{eqnarray}
and 
\begin{eqnarray}
U(D \uparrow w)_v \cup U(D)_w & \equiv & [P(D \uparrow w)_v \times S(D \uparrow w)_v] \cup [P(D)_w \times S(D)_w] \nonumber \\
& = & [P(D)_v \times (S(D)_v \setminus \{w\})] \cup [P(D)_v \times S(D)_w] \cup [P(D)_w \times S(D)_w] \nonumber
\end{eqnarray}
so upon inspection \eqref{eq:detourU} is satisfied and case 4) is done.

\begin{figure}
	\begin{tikzpicture}[->,>=stealth',shorten >=1pt]
		\node (D) at (1.5,1.75) {$D$};
		\node [draw,circle,minimum size=6mm] (v1) at (0,0) {};
		\node [draw,circle,minimum size=6mm] (v2) at (1,0) {$v$};
		\node [draw,circle,minimum size=6mm] (v3) at (1,1) {};
		\node [draw,circle,minimum size=6mm] (v4) at (2,0) {$w$};
		\node [draw,circle,minimum size=6mm] (v5) at (2,1) {};
		\node [draw,circle,minimum size=6mm] (v6) at (3,0) {};
		\foreach \from/\to in {
			v1/v2, v2/v4, v4/v6}
		\draw (\from) to (\to);
		\draw (v1) [color=white,out=-30,in=-150,looseness=1] to (v6);	% vertical alignment hack
		\draw (v2) [out=105,in=-105,looseness=1] to (v3);
		\draw (v3) [out=-75,in=75,looseness=1] to (v2);
		\draw (v4) [out=105,in=-105,looseness=1] to (v5);
		\draw (v5) [out=-75,in=75,looseness=1] to (v4);
	\end{tikzpicture}
	\quad \quad
	\begin{tikzpicture}[->,>=stealth',shorten >=1pt]
		\node (D) at (1.5,1.75) {$D \uparrow v$};
		\node [draw,circle,minimum size=6mm] (v1) at (0,0) {};
		\node [draw,circle,minimum size=6mm] (v2) at (1,0) {$v$};
		\node [draw,circle,minimum size=6mm] (v3) at (1,1) {};
		\node [draw,circle,minimum size=6mm] (v4) at (2,0) {$w$};
		\node [draw,circle,minimum size=6mm] (v5) at (2,1) {};
		\node [draw,circle,minimum size=6mm] (v6) at (3,0) {};
		\foreach \from/\to in {
			v1/v3, v3/v4, v4/v6}
		\draw (\from) to (\to);
		\draw (v1) [color=white,out=-30,in=-150,looseness=1] to (v6);	% vertical alignment hack
		\draw (v1) [out=-30,in=-150,looseness=1] to (v4);
		\draw (v4) [out=105,in=-105,looseness=1] to (v5);
		\draw (v5) [out=-75,in=75,looseness=1] to (v4);
	\end{tikzpicture}
	\quad \quad
	\begin{tikzpicture}[->,>=stealth',shorten >=1pt]
		\node (D) at (1.5,1.75) {$D \uparrow w$};
		\node [draw,circle,minimum size=6mm] (v1) at (0,0) {};
		\node [draw,circle,minimum size=6mm] (v2) at (1,0) {$v$};
		\node [draw,circle,minimum size=6mm] (v3) at (1,1) {};
		\node [draw,circle,minimum size=6mm] (v4) at (2,0) {$w$};
		\node [draw,circle,minimum size=6mm] (v5) at (2,1) {};
		\node [draw,circle,minimum size=6mm] (v6) at (3,0) {};
		\foreach \from/\to in {
			v1/v2, v2/v5, v5/v6}
		\draw (\from) to (\to);
		\draw (v1) [color=white,out=-30,in=-150,looseness=1] to (v6);	% vertical alignment hack
		\draw (v2) [out=-30,in=-150,looseness=1] to (v6);
		\draw (v2) [out=105,in=-105,looseness=1] to (v3);
		\draw (v3) [out=-75,in=75,looseness=1] to (v2);
	\end{tikzpicture}
	\quad \quad
	\begin{tikzpicture}[->,>=stealth',shorten >=1pt]
		\node (D) at (1.5,1.75) {$(D \uparrow v) \uparrow w = (D \uparrow w) \uparrow v$};
		\node [draw,circle,minimum size=6mm] (v1) at (0,0) {};
		\node [draw,circle,minimum size=6mm] (v2) at (1,0) {$v$};
		\node [draw,circle,minimum size=6mm] (v3) at (1,1) {};
		\node [draw,circle,minimum size=6mm] (v4) at (2,0) {$w$};
		\node [draw,circle,minimum size=6mm] (v5) at (2,1) {};
		\node [draw,circle,minimum size=6mm] (v6) at (3,0) {};
		\foreach \from/\to in {
			v1/v3, v1/v5, v3/v6, v5/v6}
		\draw (\from) to (\to);
		\draw (v1) [out=-30,in=-150,looseness=1] to (v6);
	\end{tikzpicture}
	\caption{\label{fig:Case4Cartoon} Cartoon for case 4 of the lemma.} 
\end{figure}

\emph{Case 5: $w \in V(D)_v^0$.}  In this case we have the following identities:
\begin{eqnarray}
P(D \uparrow v)_w & = & P(D)_w; \nonumber \\
S(D \uparrow v)_w & = & S(D)_w; \nonumber \\
P(D \uparrow w)_v & = & P(D)_v; \nonumber \\
S(D \uparrow w)_v & = & S(D)_v. \nonumber
\end{eqnarray}
From these \eqref{eq:detourU} follows trivially, so case 5) is done. $\Box$

\section{\label{sec:CommutativeContractions}Case analysis proof of $(D \uparrow u) / \{v,w\} = (D / \{v,w\}) \uparrow u$ for digraphs}

The result is trivial unless $v$ and $w$ belong to different sets of the form $V(D)_u^\bullet$. It also suffices to show the result for a modified contraction operation (denoted $\merge$ below) that yields identical copies of contracted vertices (note that this is essentially the same technical simplifcation as dealing with detours instead of bypasses). By symmetry, we need only consider the six cases where $(v,w)$ is an element of one of the following products: $V(D)_u^- \times V(D)_u^\pm$, $V(D)_u^- \times V(D)_u^+$, $V(D)_u^- \times V(D)_u^0$, $V(D)_u^\pm \times V(D)_u^+$, $V(D)_u^\pm \times V(D)_u^0$, and $V(D)_u^+ \times V(D)_u^0$.

Consider for instance the first of these cases, where $v \in V(D)_u^-$ and $w \in V(D)_u^\pm$. Writing $vw$ and $vw'$ for the identical copies of contracted vertices, $V_{u;v,w}^\bullet$ as a temporary shorthand for $V(D)_u^\bullet \setminus \{v,w\}$, and e.g., $D \merge \{v,w\}$ for the modified contraction, we have the following adjacency matrices (with irrelevant entries omitted):
\begin{equation}
\bordermatrix{\mu_D & u & v & w & V_{u;v,w}^- & V_{u;v,w}^\pm & V_{u;v,w}^+ & V_{u;v,w}^0 \cr
	u & 0 & 0 & 1 & 0 & 1 & 1 & 0 \cr
	v & 1 & 0 & \cdot & \mu_{v;-} & \cdot & \cdot & \mu_{v;0} \cr
	w & 1 & \cdot & 0 & \mu_{w;-} & \cdot & \cdot & \mu_{w;0} \cr
	V_{u;v,w}^- & 1 & \cdot & \cdot & \cdot & \cdot & \cdot & \cdot \cr
	V_{u;v,w}^\pm & 1 & \cdot & \cdot & \cdot & \cdot & \cdot & \cdot \cr
	V_{u;v,w}^+ & 0 & \mu_{+;v} & \mu_{+;w} & \cdot & \cdot & \cdot & \cdot \cr
	V_{u;v,w}^0 & 0 & \mu_{0;v} & \mu_{0;w} & \cdot & \cdot & \cdot & \cdot \cr
	} \nonumber
\end{equation}
\begin{equation}
\bordermatrix{\mu_{D \uparrow u} & u & v & w & V_{u;v,w}^- & V_{u;v,w}^\pm & V_{u;v,w}^+ & V_{u;v,w}^0 \cr
	u & 0 & 0 & 0 & 0 & 0 & 0 & 0 \cr
	v & 0 & 0 & 1 & \mu_{v;-} & 1 & 1 & \mu_{v;0} \cr
	w & 0 & \cdot & 0 & \mu_{w;-} & 1 & 1 & \mu_{w;0} \cr
	V_{u;v,w}^- & 0 & \cdot & 1 & \cdot & 1 & 1 & \cdot \cr
	V_{u;v,w}^\pm & 0 & \cdot & 1 & \cdot & 1-I & 1 & \cdot \cr
	V_{u;v,w}^+ & 0 & \mu_{+;v} & \mu_{+;w} & \cdot & \cdot & \cdot & \cdot \cr
	V_{u;v,w}^0 & 0 & \mu_{0;v} & \mu_{0;w} & \cdot & \cdot & \cdot & \cdot \cr
	} \nonumber
\end{equation}
\begin{equation}
\bordermatrix{\mu_{D \merge \{v,w\}} & u & vw & vw' & V_{u;v,w}^- & V_{u;v,w}^\pm & V_{u;v,w}^+ & V_{u;v,w}^0 \cr
	u & 0 & 1 & 1 & 0 & 1 & 1 & 0 \cr
	vw & 1 & 0 & 0 & \mu_{v;-} \lor \mu_{w;-} & \cdot & \cdot & \mu_{v;0} \lor \mu_{w;0} \cr
	vw' & 1 & 0 & 0 & \mu_{v;-} \lor \mu_{w;-} & \cdot & \cdot & \mu_{v;0} \lor \mu_{w;0} \cr
	V_{u;v,w}^- & 1 & \cdot & \cdot & \cdot & \cdot & \cdot & \cdot \cr
	V_{u;v,w}^\pm & 1 & \cdot & \cdot & \cdot & \cdot & \cdot & \cdot \cr
	V_{u;v,w}^+ & 0 & \mu_{+;v} \lor \mu_{+;w} & \mu_{+;v} \lor \mu_{+;w} & \cdot & \cdot & \cdot & \cdot \cr
	V_{u;v,w}^0 & 0 & \mu_{0;v} \lor \mu_{0;w} & \mu_{0;v} \lor \mu_{0;w} & \cdot & \cdot & \cdot & \cdot \cr
	} \nonumber
\end{equation}
whereupon both $(D \uparrow u) / \{v,w\}$ and $(D/\{v,w\}) \uparrow u$ can be seen to have the adjacency matrix
\begin{equation}
\bordermatrix{\ & u & vw & vw' & V_{u;v,w}^- & V_{u;v,w}^\pm & V_{u;v,w}^+ & V_{u;v,w}^0 \cr
	u & 0 & 0 & 0 & 0 & 0 & 0 & 0 \cr
	vw & 0 & 0 & 0 & \mu_{v;-} \lor \mu_{w;-} & 1 & 1 & \mu_{v;0} \lor \mu_{w;0} \cr
	vw' & 0 & 0 & 0 & \mu_{v;-} \lor \mu_{w;-} & 1 & 1 & \mu_{v;0} \lor \mu_{w;0} \cr
	V_{u;v,w}^- & 0 & 1 & 1 & \cdot & 1 & 1 & \cdot \cr
	V_{u;v,w}^\pm & 0 & 1 & 1 & \cdot & 1-I & 1 & \cdot \cr
	V_{u;v,w}^+ & 0 & \mu_{+;v} \lor \mu_{+;w} & \mu_{+;v} \lor \mu_{+;w} & \cdot & \cdot & \cdot & \cdot \cr
	V_{u;v,w}^0 & 0 & \mu_{0;v} \lor \mu_{0;w} & \mu_{0;v} \lor \mu_{0;w} & \cdot & \cdot & \cdot & \cdot \cr
	} \nonumber
\end{equation}
and this case is done. The other cases are entirely similar (in fact, the first, second, fourth and fifth cases are nearly identical). $\Box$

\section{\label{sec:Renormalization}Remark on renormalization}

We recall two theorems described in \cite{FriezeKaronski} regarding $D_{n,p}$:

\

\textsc{Theorem.} If $c >1$ is constant, then with high probability $D_{n,c/n}$ contains a unique strong component of size $\approx (1-\frac{x}{c})^2 n$, where $x < 1$ solves $xe^{-x} = ce^{-c}$. Furthermore, all other strong components are of logarithmic size. $\Box$

\

\textsc{Theorem.} $\lim_n \mathbb{P}(D_{n,p} \text{ is strongly connected}) = \exp(-2e^{-\lim_n (pn-\log n)})$. $\Box$

\

\begin{figure}[htbp]
\includegraphics[trim = 15mm 0mm 15mm 0mm, clip, width=180mm,keepaspectratio]{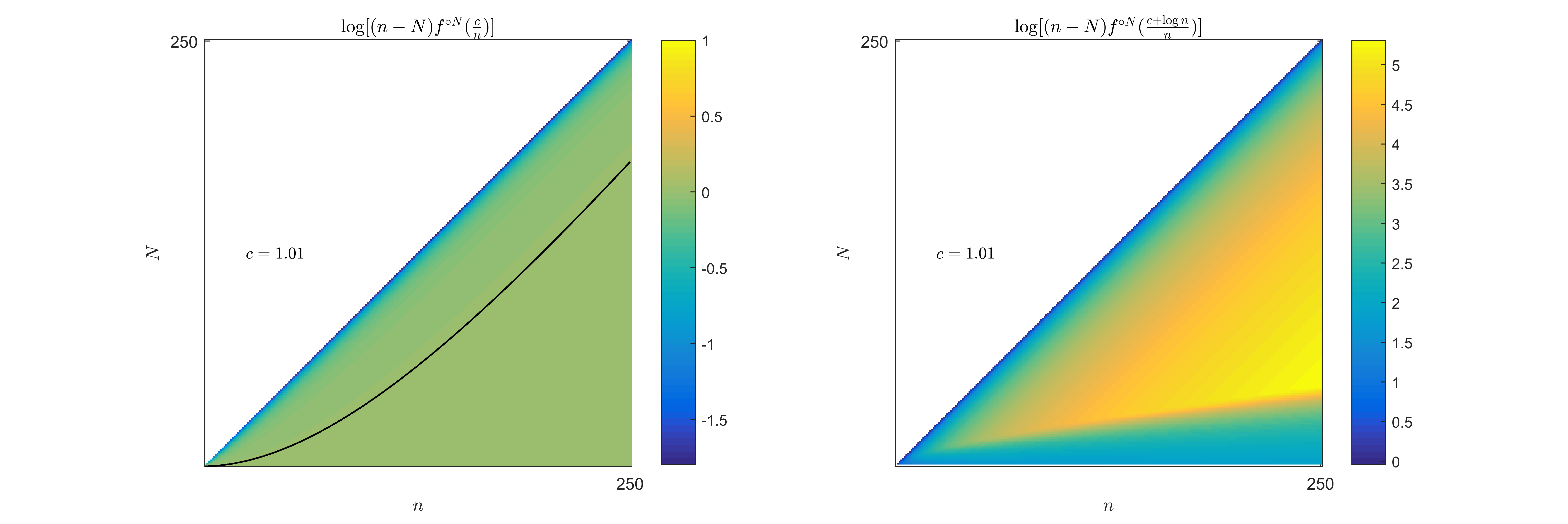}% Here is how to import pix
\caption{ \label{fig:c101} (L) $\log [(n-N) f^{\circ N}(\frac{c}{n})]$ for $c = 1.01$, with contour at 0 drawn. (R) $\log [(n-N) f^{\circ N}(\frac{c+\log n}{n})]$ for $c = 1.01$.} 
\end{figure}

\begin{figure}[htbp]
\includegraphics[trim = 15mm 0mm 15mm 0mm, clip, width=180mm,keepaspectratio]{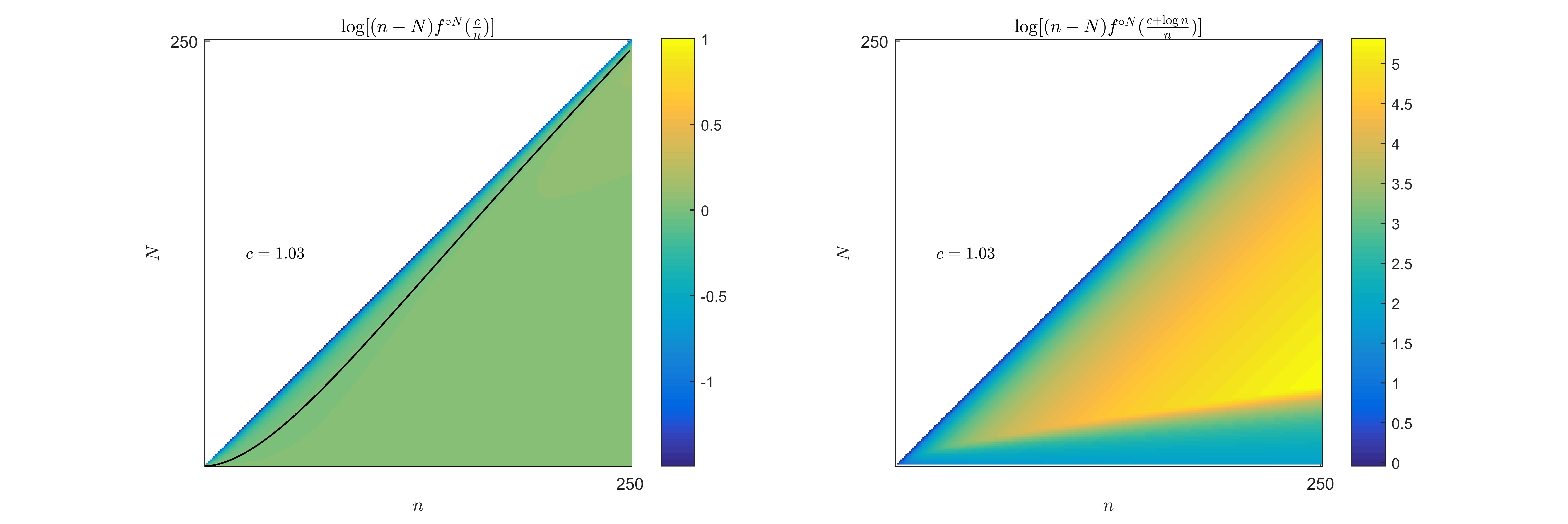}% Here is how to import pix
\caption{ \label{fig:c103} (L) $\log [(n-N) f^{\circ N}(\frac{c}{n})]$ for $c = 1.03$, with contour at 0 drawn. (R) $\log [(n-N) f^{\circ N}(\frac{c+\log n}{n})]$ for $c = 1.03$.} 
\end{figure}

These theorems suggest studying the behavior of $(n-N) f^{\circ N}(\frac{c}{n})$ and $(n-N) f^{\circ N}(\frac{c+\log n}{n})$ for $c$ constant and $0 \le N < n$. Numerics indicate that for $c > 1$ the first of these is greater than unity except for $N \approx n$, and the second is always greater than unity: see figures \ref{fig:c101} and \ref{fig:c103}.

    % N = 250;
    % c = 1.01;
    % test1 = nan(N);
    % test2 = nan(N);
    % for n = 1:N
    %     k = 0;
    %     p1 = c/n;
    %     p2 = (c+log(n))/n;
    %     f1 = p1;    % k = 0
    %     f2 = p2;
    %     test1(k+1,n) = (n-k)*f1;
    %     for k = 1:(n-1)
    %         f1 = f1^2+(1-f1^2)*f1;
    %         f2 = f2^2+(1-f2^2)*f2;
    %         test1(k+1,n) = (n-k)*f1;
    %         test2(k+1,n) = (n-k)*f2;
    %     end
    % end
    % %%
    % % redblue = [linspace(0,1,64);zeros(1,64);linspace(1,0,64)]';
    % figure;
    % subplot(1,2,1); 
    % hold on; 
    % pcolorfull(log(test1)); 
    % contour(log(test1),[0,0],'k','LineWidth',1); 
    % % colormap(redblue);
    % axis square;
    % box on;
    % text(N/10,N/2,['$c = ',num2str(c),'$'],'Interpreter','latex');
    % xlabel('$n$','Interpreter','latex');
    % ylabel('$N$','Interpreter','latex');
    % title('$\log [(n-N) f^{\circ N}(\frac{c}{n})]$','Interpreter','latex');
    % set(gca,'XTick',[0,N],'YTick',[0,N]);
    % subplot(1,2,2);
    % pcolorfull(log(test2)); 
    % % colormap(redblue);
    % axis square;
    % text(N/10,N/2,['$c = ',num2str(c),'$'],'Interpreter','latex');
    % xlabel('$n$','Interpreter','latex');
    % ylabel('$N$','Interpreter','latex');
    % title('$\log [(n-N) f^{\circ N}(\frac{c+\log n}{n})]$','Interpreter','latex');
    % set(gca,'XTick',[0,N],'YTick',[0,N]);

\end{document}